\documentclass[12pt]{article}

\usepackage[latin1]{inputenc}

\usepackage{pgf,tikz}
\usetikzlibrary{arrows}
\usetikzlibrary{patterns}

\usepackage{multirow} 
\usepackage[caption=false]{subfig}

\usepackage{amssymb}

\usepackage{amsmath,epsfig}%pstricks,pst-node}
\def\X{\textrm{X}}
\def\Y{\textrm{Y}}

\newtheorem{Theorem}{Theorem}{\bfseries}{\itshape}
{\bfseries}{\itshape}
\newtheorem{Proposition}[Theorem]{Proposition}{\bfseries}{\itshape}
{\bfseries}{\itshape}
\newtheorem{Definition}[Theorem]{Definition}{\bfseries}{\itshape}
{\bfseries}{\itshape}
{\bfseries}{\itshape}

\newcommand{\Z}{{\mathbb Z}}

\newcommand{\beq}{\begin{equation}}
\newcommand{\eeq}{\end{equation}}

\def\qed{\hfill$\Box$\\ \medskip}

\usepackage{multirow}

\def\mybox{\hbox to 12.0pt}

\def\mybigbox{\hbox to 35.0pt}
\def\myverybigbox{\hbox to 60.0pt}

\def\ov{\overline}
\def\sc{\scriptstyle}

\def\s{{\bf{s}}}

\def\wgt{{\rm wgt}}

\def\Mol{{\rm Mol}}
\def\Bac{{\rm Bac}}

\def\={\!=\!}
\def\+{\!+\!}
\def\-{\!-\!}

\def\x{{\bf x}}
\def\y{{\bf y}}
\def\z{{\bf z}}

\def\a{{\bf a}}

\def\0{{\bf 0}}

\def\1{{\bf 1}}

\def\red{\textcolor{red}}
\def\blue{\textcolor{blue}}
\def\magenta{\textcolor{magenta}}

\def\cyan{\textcolor{cyan}}
\def\brown{\textcolor{brown!80!black}}
\def\green{\textcolor{green!70!black}}

\title{Factorial supersymmetric skew Schur functions and ninth variation determinantal identities}

\author{
Ang\`ele M. Foley\thanks{ 
Department of Physics and Computer Science,
Wilfrid Laurier University,
Waterloo, Ontario, N2L 3C5, Canada ({\tt ahamel@wlu.ca})}
\and 
Ronald C. King\thanks{
Mathematical Sciences, University of Southampton, 
Southampton SO17 1BJ, England ({\tt r.c.king@soton.ac.uk})}}

\begin{document}

\maketitle

\begin{abstract}
The determinantal identities of Hamel and Goulden have recently been shown to apply to
a tableau-based ninth variation of skew Schur functions.  Here we extend this approach and its results to the analogous tableau-based ninth variation of supersymmetric skew Schur functions. These tableaux are built on entries taken from an alphabet of unprimed and primed numbers and that may be ordered in a myriad of different ways, each  leading to a determinantal identity.
At the level of the ninth variation the corresponding determinantal identities are all distinct
but the original notion of supersymmetry is lost. It is shown that this can be remedied at the
level of the sixth variation involving a doubly infinite sequence of factorial parameters. Moreover 
it is shown that the resulting factorial supersymmetric skew Schur functions are independent of the
ordering of the unprimed and primed entries in the alphabet. 
\end{abstract}

\section{Introduction}\label{Sec-introduction}

Supersymmetric Schur functions trace their genesis to the 1987 paper of Berele and Regev \cite{BR87} who defined them---under the name of ($k,\ell$)-hook Schur functions---using hook Young tableaux.  Their combinatorics has been well-developed since, e.g.\ \cite{Kwon,Pragacz,Rem1,Rem2,YR}, and they have been generalized, particularly in the quasisymmetric direction \cite{MN} and the  factorial direction \cite{Mol}.  Here we further extend the factorial direction, defining a tableau-based ninth variation of factorial supersymmetric Schur functions,
along the lines of that proposed by Bachmann \cite{Bac} in a non-supersymmetric context.

Let $\x=(x_1,x_2,\ldots,x_m)$ and $\y=(y_1,y_2,\ldots,y_n)$ be two sequences of indeterminate parameters.
A function $f(\x,\y)$ is said to be supersymmetric if it is invariant under independent permutations of the 
components of $\x$ and of $\y$, and if it is independent of $t$ if one sets $x_1=t$ and $y_1=-t$~\cite{Ste85}.

For partitions $\lambda$ and $\mu$ the supersymmetric skew Schur functions $s_{\lambda/\mu}(\x,\y)$ are a 
generalisation of the symmetric skew Schur functions $s_{\lambda/\mu}(\x)$ in which the
underlying alphabet $\x$ is replaced by $\z=(\x,\y)$. Just as the symmetric nature of $s_{\lambda/\mu}(\x)$ follows
rather simply from its definition in terms of semistandard skew tableaux with entries $k$ from the 
alphabet $1<2<\cdots<m$ carrying weight $x_k$, so the supersymmetric nature of $s_{\lambda/\mu}(\x,\y)$ follows  from its 
definition in terms of super skew tableaux~\cite{BR87} with entries $k$ and $\ell'$ from the alphabet $1<2<\ldots<m<1'<2'<\cdots<n'$,
carrying weights $x_k$ and $y_{\ell}$, respectively.

Factorial skew Schur functions $s_{\lambda/\mu}(\x|\a)$ are referred to by Macdonald in~\cite{Mac92} (see also \cite{Mac95}) as his sixth variation. 
They involve a doubly infinite sequence of factorial parameters $\a=(\ldots,a_{-2},a_{-1},a_0,a_1,a_2,\ldots)$.
In this factorial case, an entry $k$ in a box on the $c$th diagonal of a semistandard skew tableau carries weight $x_k+a_c$.
Factorial supersymmetric skew Schur functions $s_{\lambda/\mu}(\x,\y|\a)$ were introduced by Molev~\cite{Mol},
with skew supertableaux entries $k$ and $\ell'$ from the alphabet $n'<\cdots<2'<1'<1<2<\cdots<m$ carrying weights $x_k+a_{k+c}$ and 
$y_\ell-a_{\ell+c}$, respectively, if they appear on the $c$th diagonal.

Introducing parameters $\X=(x_{kc})$ allows one to give a combinatorial realisation of Macdonald's ninth variation skew Schur functions 
which were originally defined algebraically. These ninth variation skew Schur functions, $s_{\lambda/\mu}(\X)$, are again
defined in terms of semistandard skew tableaux~\cite{BC,FK20}, with an entry $k$ on the $c$th diagonal carrying 
weight $x_{kc}$.
This notion has been extended to another ninth variation of Schur functions in~\cite{Bac}. This was based on certain ordered Young
tableaux that might be viewed as precursors of supertableaux. The entries $k$ on the $c$th diagonal were given a weight $F=(f_{kc})$ augmented
by factors $t$ or $(1-t)$ depending on the nature of their immediate neighbours. In the same spirit, setting
$\X=(x_{kc})$ and $\Y=(y_{\ell c})$ allows one to define ninth variation super skew Schur functions $s_{\lambda/\mu}(\X,\Y)$ 
that are based on skew supertableaux with entries $k$ and $\ell'$ on the $c$th diagonal carrying weights $x_{kc}$ and $y_{\ell c}$, respectively.
It will be shown that a special case of this yields the variation based on ordered Young tableaux.

Through the use of a non-intersecting path model of semistandard skew tableaux, Hamel and Goulden~\cite{HG}
derived a family of determinantal expansions of the skew Schur functions $s_{\lambda/\mu}(\x)$; 
one for each so-called outside decomposition. Their result is embodied in the following Theorem for which 
the full notation will be explained later.

\begin{Theorem} [Hamel and Goulden~\cite{HG}] \label{The-HG} 
Let $\lambda$ and $\mu$ be partitions such that $\mu\subseteq\lambda$, and let 
$\theta=(\theta_1,\theta_2,\ldots,\theta_s)$ be an outside decomposition of $F^{\lambda/\mu}$.
Then for all $\x=(x_1,x_2,\ldots,x_n)$
\begin{equation}\label{eqn-HGsfn}
    s_{\lambda/\mu}(\x) = \left|\,  s_{(\theta_p\#\theta_q)}(\x) \,\right|_{1\leq p,q \leq s}\,,
\end{equation}
where $(\theta_p\#\theta_q)$ is either undefined, or a single edge, or a strip
formed by adjoining the strips $\theta_p$ and $\theta_q$, with some overlap if necessary,
while preserving their shape and their content, that is to say 
the diagonals on which their boxes lie.
\end{Theorem}

It has been shown that these same identities apply to both $s_{\lambda/\mu}(\x|\a)$~\cite{Sch} 
and $s_{\lambda/\mu}(\X)$~\cite{BC,FK20}. Here we show that they also apply to the supersymmetric skew Schur functions $s_{\lambda/\mu}(\x,\y)$ 
and $s_{\lambda/\mu}(\x,\y|\a)$, as well as to the ninth variation super skew Schur functions $s_{\lambda/\mu}(\X,\Y)$.

The paper proceeds as follows.  In Section \ref{Sec-ssfmn} we define the tableaux, weighting schemes and resulting symmetric functions we will use in this paper. Section \ref{Sec-outside-decompositions}provides definitions and background on outside decompositions, cutting strips, lattices and other combinatorial constructs 
behind the construction of
determinantal identities. This section also contains the main result, a set of determinantal identities  
for the tableau-based ninth variation supersymmetric skew Schur functions. Section \ref{Sec:factssfn} concerns the specialisation to factorial supersymmetric skew Schur functions and shows {that they are 
independent of the ordering of the variables. In Section \ref{Sec-conclusion} we establish the precise connection with the results of Molev~\cite{Mol} and Bachmann~\cite{Bac}. 
}

\section{Supersymmetric Schur functions}~\label{Sec-ssfmn}

We begin with the fundamental notation.  Our definition of supersymmetric functions is in terms of a {\em partition}, a weakly decreasing series of positive integers.  Formally, a partition $\lambda=(\lambda_1, \lambda_2, \ldots, \lambda_l)$ is such that $\lambda_1 \geq \lambda_2 \geq \cdots \geq \lambda_l > 0$.  Each $\lambda_i$ is termed a {\em part} of the partition and thus partition $\lambda$ has $\ell$ parts or, we also say the length of $\lambda$ is $l$ and use the notation $\ell(\lambda)=l$.  In order to define the supersymmetric function we need a combinatorial object called a {\em Young diagram.}  This consists of an arrangement of boxes that is top and left justified and has $\lambda_k$ boxes in each row, $1 \leq k \leq l$.  We use the notation $F^\lambda$ for the Young diagram of the partition $\lambda$. We identify the box in row $i$, column $j$ by $(i,j)$.  Each box has a parameter called {\em content} associated to it. The content, $c$, of box $(i,j)$ is defined as $c=j-i$.  Note that this means the contents of boxes above the main diagonal will be positive and those below the main diagonal will be negative (instead of writing $c=-a$ for these boxes, we will use the notation $c=\ov{a}$).. The main diagonal itself has boxes of $0$ content.

Our definition of skew supersymmetric functions will require an additional partition. This second partition, $\mu=(\mu_1, \mu_2, \ldots, \mu_m)$ where $m \leq l$, is defined such that $\mu_i \leq \lambda_i$ for all $1\leq i \leq m$, and we say that $\mu\subseteq\lambda$.  The boxes corresponding to the $F^\mu$, the Young diagram of $\mu$, are removed from $F^{\lambda}$,  the Young diagram of $\lambda$, to create $F^{\lambda/\mu}$, the Young diagram of shape $\lambda/\mu$
The row and column labels, $i$ and $j$, of each box of $F^{\lambda/\mu}$ are inherited 
from those of the corresponding box of $F^\lambda$, as is its content $c=j-i$.

Given a skew Young diagram $F^{\lambda/\mu}$, let ${\cal T}^{\lambda/\mu}$ be the set of all skew supertableaux $T$ 
obtained by filling each box of $F^{\lambda/\mu}$ with an 
with an entry from the alphabet $1<2<\cdots<m<1'<2'<\cdots<n'$ 
in such a way that the entries are weakly increasing across rows from left to right with no two
identical unprimed entries $k$ in the same column and no two identical primed entries $\ell'$ in the 
same row. We denote the entry at position $(i,j)$ in the $i$th row and $j$th column of $T$ by $t_{ij}$.
We will define three different types of supersymmetric skew Schur function, 
each based on one of the variations specified  
by Macdonald \cite{Mac92}.

Let $\x=(x_1,x_2,\ldots,x_m)$ and $\y=(y_1,y_2,\ldots,y_n)$ be two sequences of indeterminates,
then first variation supersymmetric skew Schur function $s_{\lambda/\mu}(\x,\y)$ may be defined by
\begin{equation}\label{eqn-sfnxy}
   s_{\lambda/\mu}(\x,\y)= \sum_{T\in{\cal T}^{\lambda/\mu}} \prod_{(i,j)\in F^{\lambda/\mu}} \wgt(t_{ij}), 
\end{equation}
where $\wgt(t_{ij})=x_k$ if $t_{ij}=k$ for any $k\in M=\{1,2,\ldots,m\}$
and $\wgt(t_{ij})=y_\ell$ if $t_{ij}=\ell'$ for any $\ell\in N=\{1,2,\ldots,n\}$.

A tableau-based sixth variation of supersymmetric skew Schur functions is defined by introducing  
factorial parameters $\a=(\ldots,a_{-2},a_{-1},a_0,a_1,a_2,\ldots)$ 
and modifying the weights so as to give
\begin{equation}\label{eqn-sfnxya}
   s_{\lambda/\mu}(\x,\y|\a)=\sum_{T\in{\cal T}^{\lambda/\mu}} \prod_{(i,j)\in F^{\lambda/\mu}} \wgt(t_{ij}), 
\end{equation} 
where $\wgt(t_{ij})=x_k+a_{k+c}$ if $t_{ij}=k$ for any $k\in M$
and $\wgt(t_{ij})=y_\ell-a_{m-\ell+1+c}$ if $t_{ij}=\ell'$ for any $\ell\in N$, and $c=j-i$.
It will be shown later that, with this particular parameterisation of weights, the resulting factorial 
supersymmetric skew Schur functions are indeed supersymmetric.

More generally, a tableau-based ninth variation of supersymmetric skew Schur functions 
is defined by introducing weight parameters $\X=(x_{kc})$ and $\Y=(y_{\ell c})$ 
so as to produce
\begin{equation}\label{eqn-sfnXY}
   s_{\lambda/\mu}(\X,\Y)= \sum_{T\in{\cal T}^{\lambda/\mu}} \prod_{(i,j)\in F^{\lambda/\mu}} \wgt(t_{ij}),
\end{equation}
where now $\wgt(t_{ij})=x_{kc}$ if $t_{ij}=k$ for any $k\in M$
and $\wgt(t_{ij})=y_{\ell c}$ if $t_{ij}=\ell'$ for any $\ell \in N$, and $c=j-i$. These 
ninth variation skew Schur functions are not supersymmetric.

The weights $\wgt(t_{ij})$ appropriate to equations (\ref{eqn-sfnxy})--(\ref{eqn-sfnXY}) are summarized in the following table:

\begin{equation}
\begin{tabular}{|l|l|l|l|} \hline
Variation & $t_{ij}$ & $\wgt(t_{ij})$ & $\wgt(t_{ij})$ \\ \hline
first & k & $x_k$ &  \\
      &  $\ell'$ &    & $y_\ell$ \\ \hline
      sixth & k & $x_k+a_{k+c}$ &  \\
      &  $\ell'$ &    & $y_\ell-a_{m-\ell+1+c}$ \\ \hline
ninth & k & $x_{kc}$ &  \\
      &  $\ell'$ &    & $y_{\ell c}$ \\ \hline

\end{tabular}
\end{equation}

The above definitions are based on supertableaux $T$ with entries taken from the particular
alphabet $1<2<\cdots<m<1'<2'<\cdots<n'$. As a means of generalising them it is convenient to introduce 
alternative skew supertableaux as follows.

\begin{Definition}\label{Def-R'tableau}
Let $R=M\dot\cup N=\{1<2<\cdots<m+n\}$ be an alphabet that is the disjoint union of
sub-alphabets $M=\{i_1<i_2<\cdots<i_m\}$ and $N=\{j_1<j_2<\cdots<j_n\}$,
and let $R'=M\dot\cup N'$ be a marked version of $R$ obtained by leaving any $r\in M$ unprimed 
but adding a prime to any $r\in N$ to give $r'\in N'$. 
Given a skew Young diagram $F^{\lambda/\mu}$, let $T\in{\cal T}_{R'}^{\lambda/\mu}$ be the set of all
skew supertableaux $T$ obtained by filling each box of $F^{\lambda/\mu}$, with an
entry from $R'$ in such a way that
\begin{itemize}
\item entries are weakly increasing across rows from left to right and down columns from top to bottom;
\item no two identical unprimed entries appear in the same column;
\item no two identical primed entries appear in the same row.
\end{itemize}
\end{Definition}

For example, if $R=\{\blue1<\blue2<\red<\blue4<\red5<\red6<\red7<\blue8<\blue9\}$
with $M=\{\blue1<\blue2<\blue4<\blue8<\blue9\}$ and $N=\{\red3<\red5<\red6<\red7\}$
then $R'=\{\blue1<\blue2<\red{3'}<\blue4<\red{5'}<\red{6'}<\red{7'}<\blue8<\blue9\}$
then a typical skew supertableau $T\in{\cal T}_{R'}^{655443/32}$ takes the form
\begin{equation}\label{eqn-tab-ex}
\vcenter{\hbox{
\begin{tikzpicture}[x={(0in,-0.2in)},y={(0.2in,0in)}] % matrix coordinate
%shape lambda
\foreach \j in {1,...,6} \draw[thick] (1,\j) rectangle +(-1,-1); 
\foreach \j in {1,...,5} \draw[thick] (2,\j) rectangle +(-1,-1);
\foreach \j in {1,...,5} \draw[thick] (3,\j) rectangle +(-1,-1);
\foreach \j in {1,...,4} \draw[thick] (4,\j) rectangle +(-1,-1);
\foreach \j in {1,...,4} \draw[thick] (5,\j) rectangle +(-1,-1);
\foreach \j in {1,...,3} \draw[thick] (6,\j) rectangle +(-1,-1);
%shape mu
\foreach \j in {1,...,3} \draw(1-0.5,\j-0.5)node{$\ast$};
\foreach \j in {1,...,2} \draw(2-0.5,\j-0.5)node{$\ast$};
%content of lambda/mu
\draw(1-0.5,3+0.5) node{$\red{3'}$};\draw(1-0.5,4+0.5) node{$\red{5'}$};\draw(1-0.5,5+0.5) node{$\blue{8}$};
\draw(2-0.5,2+0.5) node{$\blue{2}$};\draw(2-0.5,3+0.5) node{$\red{3'}$};\draw(2-0.5,4+0.5) node{$\red{5'}$};
\draw(3-0.5,0+0.5) node{$\blue{1}$};\draw(3-0.5,1+0.5) node{$\blue{1}$};\draw(3-0.5,2+0.5) node{$\red{3'}$};\draw(3-0.5,3+0.5) node{$\blue{4}$};\draw(3-0.5,4+0.5) node{$\blue{9}$};
\draw(4-0.5,0+0.5) node{$\red{3'}$};\draw(4-0.5,1+0.5) node{$\blue{4}$};\draw(4-0.5,2+0.5) node{$\blue{4}$};\draw(4-0.5,3+0.5) node{$\red{6'}$};
\draw(5-0.5,0+0.5) node{$\red{3'}$};\draw(5-0.5,1+0.5) node{$\red{7'}$};\draw(5-0.5,2+0.5) node{$\blue{8}$};\draw(5-0.5,3+0.5) node{$\blue{8}$};
\draw(6-0.5,0+0.5) node{$\blue{4}$};\draw(6-0.5,1+0.5) node{$\blue{8}$};\draw(6-0.5,2+0.5) node{$\blue{9}$};
\end{tikzpicture}
}}
\end{equation}

With this definition of skew supertableau based on the alphabet $R'$, the corresponding ninth variation 
supersymmetric skew Schur functions, or generalised skew Schur functions, are defined as follows:
\begin{Definition}\label{Def-sfnXY-R'}
Let $R'=M\dot\cup N'$ be the marked version of the alphabet $R=M\dot\cup N$ with 
$M=\{i_1<i_2<\cdots<i_m\}$ and $N=\{j_1<j_2<\cdots<j_n\}$, as described in Definition~\ref{Def-R'tableau}. Then 
\begin{equation}\label{eqn-sfnXY-R'}
   s_{\lambda/\mu}^{R'}(\X,\Y)= \sum_{T\in{\cal T}_{R'}^{\lambda/\mu}} \prod_{(i,j)\in F^{\lambda/\mu}} \wgt(t_{ij}),
\end{equation}
where $\wgt(t_{ij})=x_{kc}$ if $t_{ij}=i_k\in M$
and $\wgt(t_{ij})=y_{\ell c}$ if $t_{ij}=j'_\ell\in N'$, and $c=j-i$. 
\end{Definition}

It follows from this definition that, in the case $\lambda=(6,5,5,4,4,3)$ and $\mu=(3,2)$,
with $M=\{1<2<4<8<9\}$ and $N=\{3<5<6<7\}$, the contribution 
to $s_{\lambda/\mu}^{R'}(\X,\Y)$ of the skew supertableau $T$ of (\ref{eqn-tab-ex}) is the product of the entries 
in the following array:
\begin{equation}\label{eqn-wgt-ex}
\vcenter{\hbox{
\begin{tikzpicture}[x={(0in,-0.2in)},y={(0.4in,0in)}] % matrix coordinate
%shape lambda
\foreach \j in {1,...,6} \draw[thick] (1,\j) rectangle +(-1,-1); 
\foreach \j in {1,...,5} \draw[thick] (2,\j) rectangle +(-1,-1);
\foreach \j in {1,...,5} \draw[thick] (3,\j) rectangle +(-1,-1);
\foreach \j in {1,...,4} \draw[thick] (4,\j) rectangle +(-1,-1);
\foreach \j in {1,...,4} \draw[thick] (5,\j) rectangle +(-1,-1);
\foreach \j in {1,...,3} \draw[thick] (6,\j) rectangle +(-1,-1);
%shape mu
\foreach \j in {1,...,3} \draw(1-0.5,\j-0.5)node{$\ast$};
\foreach \j in {1,...,2} \draw(2-0.5,\j-0.5)node{$\ast$};
%content of lambda/mu
\draw(1-0.5,3+0.5) node{$\red{y_{13}}$};\draw(1-0.5,4+0.5) node{$\red{y_{24}}$};\draw(1-0.5,5+0.5) node{$\blue{x_{45}}$};
\draw(2-0.5,2+0.5) node{$\blue{x_{21}}$};\draw(2-0.5,3+0.5) node{$\red{y_{12}}$};\draw(2-0.5,4+0.5) node{$\red{y_{23}}$};
\draw(3-0.5,0+0.5) node{$\blue{x_{1\ov2}}$};\draw(3-0.5,1+0.5) node{$\blue{x_{1\ov1}}$};\draw(3-0.5,2+0.5) node{$\red{y_{10}}$};\draw(3-0.5,3+0.5) node{$\blue{x_{31}}$};\draw(3-0.5,4+0.5) node{$\blue{x_{52}}$};
\draw(4-0.5,0+0.5) node{$\red{y_{1\ov3}}$};\draw(4-0.5,1+0.5) node{$\blue{x_{3\ov2}}$};\draw(4-0.5,2+0.5) node{$\blue{x_{3\ov1}}$};\draw(4-0.5,3+0.5) node{$\red{ y_{30}}$};
\draw(5-0.5,0+0.5) node{$\red{y_{1\ov4}}$};\draw(5-0.5,1+0.5) node{$\red{y_{4\ov3}}$};\draw(5-0.5,2+0.5) node{$\blue{x_{4\ov2}}$};\draw(5-0.5,3+0.5) node{$\blue{x_{4\ov1}}$};
\draw(6-0.5,0+0.5) node{$\blue{x_{3\ov5}}$};\draw(6-0.5,1+0.5) node{$\blue{x_{4\ov4}}$};\draw(6-0.5,2+0.5) node{$\blue{x_{5\ov3}}$};
\end{tikzpicture}
}}
\end{equation}
Notice that all the entries are distinct because no two identical entries, whether unprimed or primed, may appear on the same diagonal.

\section{Outside decompositions and cutting strips}~\label{Sec-outside-decompositions}

To begin we review the definitions and notation from Hamel and Goulden \cite{HG}.  Their main mechanism is an {\em outside decomposition} which is a specific partitioning of a (skew) Young diagram into smaller pieces.  To define this formally, first define a {\em strip} in a skew Young diagram as a set of edgewise connected boxes that contains no $2 \times 2$ set of boxes.  One artifact of this definition is that each diagonal is represented at most once in any given strip.  Strips have a direction to them:  we consider that they start at the box of lowest content, i.e. their leftmost and bottommost box, which we term that the {\em starting box} of the strip, and that they terminate at the box of highest content, i.e., their rightmost and topmost box, which we term that the {\em ending box} of the strip.

Now we consider partitioning the skew Young diagram into strips such that each box belongs to exactly one strip. We add an additional restriction that the starting box of every strip must be on the left or bottom boundary of the skew Young diagram while the ending box of every strip must be on the right or top boundary of the skew Young diagram (note that because the diagram is skew this of course includes the boundary in the upper left where the Young diagram for $\mu$ was removed).  The disjoint union $\Theta$ of such a set of strips $(\theta_1,\theta_2,\ldots,\theta_s)$ was named an {\em outside decomposition} by Hamel and Goulden \cite{HG}.  Note that we consider the content of a strip in an outside decomposition to be the content the strip had in the skew Young diagram it came from.  Given a skew Young diagram there are potentially many possible outside decompositions.

Additionally, Chen et al. \cite{CYY}, introduced a way of looking at an outside decomposition from the perspective of a type of ``master'' strip which determines the shape of all the other strips.  They called this master the {\em cutting strip}.  One pictures the cutting strip as moving through the skew Young diagram ``cutting'' it at successive steps and leaving the outside decomposition strips in its wake.  More formally, the cutting strip, $\phi$, is a strip with one box of each content number that appears in the skew Young diagram $F^{\lambda/\mu}$.  To cut $F^{\lambda/\mu}$ into disjoint strips for an outside decomposition superimpose diagonally translated copies of the cutting strip onto $F^{\lambda/\mu}$ and collect the resulting strips.  An arbitrary strip $\theta_p$  with  starting box of content $a$ and  ending box of content $b$, where $a  \leq b+1$, and created in this way from a cutting strip, shares some characteristics with $\phi$: it will have the same shape and content as what we will call $\phi_{ab}$,  the portion of $\phi$ from the box of content $a$ to the box of content $b$ (recall that $\phi$ includes boxes of all contents in $F^{\lambda/\mu}$ so this identification is always possible).   We write $\theta_p\simeq\phi_{ab}$.

If the cutting strip $\phi$ itself extends from a box of content $a$ to one of content $b$ then it is composed of boxes 
of content $a,a+1,\ldots,b$ and these values may be be divided into two sets $L$ and $B$ according as the box of content $c-1$ 
lies to the left or below the box of content $c$. The value $c=a$ may be assigned to $L$ or $B$ as is most convenient.
In the case $\lambda=(6,5,5,4,4,3)$ and $\mu=(3,2)$ this is all illustrated below
for one particular choice of cutting strip $\phi$:
\begin{equation}\label{eqn-cs}
\begin{array}{cc}
\mbox{Content of $F^{\lambda/\mu}$}&\mbox{Content of cutting strip $\phi$}\cr
%\cr
\vcenter{\hbox{
\begin{tikzpicture}[x={(0in,-0.2in)},y={(0.2in,0in)}] % matrix coordinate
%shape lambda
\foreach \j in {1,...,6} \draw[thick] (1,\j) rectangle +(-1,-1); 
\foreach \j in {1,...,5} \draw[thick] (2,\j) rectangle +(-1,-1);
\foreach \j in {1,...,5} \draw[thick] (3,\j) rectangle +(-1,-1);
\foreach \j in {1,...,4} \draw[thick] (4,\j) rectangle +(-1,-1);
\foreach \j in {1,...,4} \draw[thick] (5,\j) rectangle +(-1,-1);
\foreach \j in {1,...,3} \draw[thick] (6,\j) rectangle +(-1,-1);
%shape mu
\foreach \j in {1,...,3} \draw(1-0.5,\j-0.5)node{$\ast$};
\foreach \j in {1,...,2} \draw(2-0.5,\j-0.5)node{$\ast$};
%content of lambda/mu
\foreach \k in {3,...,5} \draw(1-0.5,\k+0.5) node{$\k$};
\foreach \k in {1,...,3} \draw(2-0.5,\k+1.5) node{$\k$};
\foreach \k in {1,...,2} \draw(3-0.5,-\k+2.5) node{$\ov{\k}$};\foreach \k in {0,...,2} \draw(3-0.5,\k+2.5) node{$\k$};
\foreach \k in {1,...,3} \draw(4-0.5,-\k+3.5) node{$\ov{\k}$};\foreach \k in {0,...,0} \draw(4-0.5,\k+3.5) node{$\k$};
\foreach \k in {1,...,4} \draw(5-0.5,-\k+4.5) node{$\ov{\k}$};
\foreach \k in {3,...,5} \draw(6-0.5,-\k+5.5) node{$\ov{\k}$};
\end{tikzpicture}
}}
&%%cutting strip
\vcenter{\hbox{
\begin{tikzpicture}[x={(0in,-0.2in)},y={(0.2in,0in)}] % matrix coordinate
%shape phi
\foreach \j in {6,...,7} \draw[thick] (1,\j) rectangle +(-1,-1);
\foreach \j in {3,...,6} \draw[thick] (2,\j) rectangle +(-1,-1);
\foreach \j in {2,...,3} \draw[thick] (3,\j) rectangle +(-1,-1);
\foreach \j in {2,...,2} \draw[thick] (4,\j) rectangle +(-1,-1);
\foreach \j in {1,...,2} \draw[thick] (5,\j) rectangle +(-1,-1);
%content of phi
\draw(1-0.5,6-0.5) node {$\green{4}$}; \draw(1-0.5,7-0.5) node {$\brown{5}$};  %\draw(1-0.5,8-0.5) node {$\brown{6}$};
\draw(2-0.5,3-0.5) node {$\green{0}$}; \draw(2-0.5,4-0.5) node {$\brown{1}$}; \draw(2-0.5,5-0.5) node {$\brown{2}$}; \draw(2-0.5,6-0.5) node {$\brown{3}$};
\draw(3-0.5,2-0.5) node {$\green{\ov2}$}; \draw(3-0.5,3-0.5) node {$\brown{\ov1}$};
\draw(4-0.5,2-0.5) node {$\green{\ov3}$};
\draw(5-0.5,1-0.5) node {$\green{\ov5}$};\draw(5-0.5,2-0.5) node {$\brown{\ov4}$};
\end{tikzpicture}
}}
\cr\cr
\mbox{Superposition of cutting strips}&\mbox{Outside decomposition}\cr\cr
\vcenter
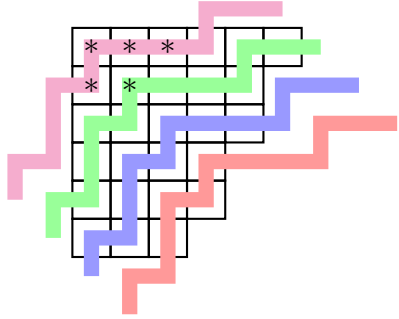
{\hbox{
\begin{tikzpicture}[x={(0in,-0.2in)},y={(0.2in,0in)}] % matrix coordinate
%shape lambda
\foreach \j in {1,...,6} \draw[thick] (1,\j) rectangle +(-1,-1); 
\foreach \j in {1,...,5} \draw[thick] (2,\j) rectangle +(-1,-1);
\foreach \j in {1,...,5} \draw[thick] (3,\j) rectangle +(-1,-1);
\foreach \j in {1,...,4} \draw[thick] (4,\j) rectangle +(-1,-1);
\foreach \j in {1,...,4} \draw[thick] (5,\j) rectangle +(-1,-1);
\foreach \j in {1,...,3} \draw[thick] (6,\j) rectangle +(-1,-1);
%cutting strips
\fill[blue!40!white](6.5+0,0.3+0)--(5.3+0,0.3+0)--(5.3+0,1.3+0)--(3.3+0,1.3+0)--(3.3+0,2.3+0)--(2.3+0,2.3+0)--(2.3+0,5.3+0)--(1.3+0,5.3+0)--(1.3+0,7.5+0)
--(1.7+0,7.5+0)--(1.7+0,5.7+0)--(2.7+0,5.7+0)--(2.7+0,2.7+0)--(3.7+0,2.7+0)--(3.7+0,1.7+0)--(5.7+0,1.7+0)--(5.7+0,0.7+0)--(6.5+0,0.7+0)--cycle;
\fill[red!40!white](6.5+1,0.3+1)--(5.3+1,0.3+1)--(5.3+1,1.3+1)--(3.3+1,1.3+1)--(3.3+1,2.3+1)--(2.3+1,2.3+1)--(2.3+1,5.3+1)--(1.3+1,5.3+1)--(1.3+1,7.5+1)
--(1.7+1,7.5+1)--(1.7+1,5.7+1)--(2.7+1,5.7+1)--(2.7+1,2.7+1)--(3.7+1,2.7+1)--(3.7+1,1.7+1)--(5.7+1,1.7+1)--(5.7+1,0.7+1)--(6.5+1,0.7+1)--cycle;
\fill[green!40!white](6.5-1,0.3-1)--(5.3-1,0.3-1)--(5.3-1,1.3-1)--(3.3-1,1.3-1)--(3.3-1,2.3-1)--(2.3-1,2.3-1)--(2.3-1,5.3-1)--(1.3-1,5.3-1)--(1.3-1,7.5-1)
--(1.7-1,7.5-1)--(1.7-1,5.7-1)--(2.7-1,5.7-1)--(2.7-1,2.7-1)--(3.7-1,2.7-1)--(3.7-1,1.7-1)--(5.7-1,1.7-1)--(5.7-1,0.7-1)--(6.5-1,0.7-1)--cycle;
\fill[magenta!40!white](6.5-2,0.3-2)--(5.3-2,0.3-2)--(5.3-2,1.3-2)--(3.3-2,1.3-2)--(3.3-2,2.3-2)--(2.3-2,2.3-2)--(2.3-2,5.3-2)--(1.3-2,5.3-2)--(1.3-2,7.5-2)
--(1.7-2,7.5-2)--(1.7-2,5.7-2)--(2.7-2,5.7-2)--(2.7-2,2.7-2)--(3.7-2,2.7-2)--(3.7-2,1.7-2)--(5.7-2,1.7-2)--(5.7-2,0.7-2)--(6.5-2,0.7-2)--cycle;
%shape mu
\foreach \j in {1,...,3} \draw(1-0.5,\j-0.5)node{$\ast$};
\foreach \j in {1,...,2} \draw(2-0.5,\j-0.5)node{$\ast$};
\end{tikzpicture}
}}
&%outside decomposition
\vcenter
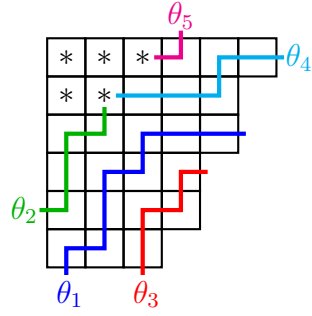
{\hbox{
\begin{tikzpicture}[x={(0in,-0.2in)},y={(0.2in,0in)}] % matrix coordinate
%shape lambda
\foreach \j in {1,...,6} \draw[thick] (1,\j) rectangle +(-1,-1); 
\foreach \j in {1,...,5} \draw[thick] (2,\j) rectangle +(-1,-1);
\foreach \j in {1,...,5} \draw[thick] (3,\j) rectangle +(-1,-1);
\foreach \j in {1,...,4} \draw[thick] (4,\j) rectangle +(-1,-1);
\foreach \j in {1,...,4} \draw[thick] (5,\j) rectangle +(-1,-1);
\foreach \j in {1,...,3} \draw[thick] (6,\j) rectangle +(-1,-1);
%shape mu
\foreach \j in {1,...,3} \draw(1-0.5,\j-0.5)node{$\ast$};
\foreach \j in {1,...,2} \draw(2-0.5,\j-0.5)node{$\ast$};
% Theta labels  
\draw(6.7,0.6)node{$\blue{\theta_1}$};
\draw(4.5,-0.6)node{$\green{\theta_2}$};
\draw(6.7,2.6)node{$\red{\theta_3}$};
\draw(0.5,6.6)node{$\cyan{\theta_4}$};
\draw(-0.6,3.5)node{$\magenta{\theta_5}$};
%\theta 
\draw[blue,ultra thick](6.2,0.5)--(5.5,0.5)--(5.5,1.5)--(3.5,1.5)--(3.5,2.5)--(2.5,2.5)--(2.5,5.2);%--(1.5,5.5)--(1.5,7.2);
\draw[green!70!black,ultra thick](4.5,-0.2)--(4.5,0.5)--(2.5,0.5)--(2.5,1.5)--(1.8,1.5);
\draw[red,ultra thick](6.2,2.5)--(4.5,2.5)--(4.5,3.5)--(3.5,3.5)--(3.5,4.2);
\draw[cyan,ultra thick](1.5,1.8)--(1.5,4.5)--(0.5,4.5)--(0.5,6.2);
\draw[magenta,ultra thick](0.5,2.8)--(0.5,3.5)--(-0.2,3.5);
\end{tikzpicture}
}}
\end{array}
\end{equation}
In the above example we have distinguished those boxes of $\phi$ having content in $L$ or $B$, by rendering
these values in brown or green, respectively. We thus have $\phi=\phi_{\ov55}$ with $L=\{\brown{\ov4,\ov1,1,2,3,5}\}$ and 
$B=\{\green{\ov5,\ov3,\ov2,0,4}\}$.

It follows in this example that the strips $\theta_p$ and their content are given, along with their identification 
with portions of $\phi$, by:
\begin{equation}\label{eqn-strips}
\begin{array}{c}
\theta_1=
\vcenter{\hbox{
\begin{tikzpicture}[x={(0in,-0.2in)},y={(0.2in,0in)}] % matrix coordinate
%shape phi
%\foreach \j in {6,...,7} \draw[thick] (1,\j) rectangle +(-1,-1);
\foreach \j in {3,...,5} \draw[thick] (2,\j) rectangle +(-1,-1);
\foreach \j in {2,...,3} \draw[thick] (3,\j) rectangle +(-1,-1);
\foreach \j in {2,...,2} \draw[thick] (4,\j) rectangle +(-1,-1);
\foreach \j in {1,...,2} \draw[thick] (5,\j) rectangle +(-1,-1);
%content of phi
\draw(5-0.5,1-0.5) node {$\blue{\ov5}$};
\draw(5-0.5,2-0.5) node {$\blue{\ov4}$};
\draw(4-0.5,2-0.5) node {$\blue{\ov3}$};
\draw(3-0.5,2-0.5) node {$\blue{\ov2}$}; 
\draw(3-0.5,3-0.5) node {$\blue{\ov1}$};
\draw(2-0.5,3-0.5) node {$\blue{0}$}; 
\draw(2-0.5,4-0.5) node {$\blue{1}$}; 
\draw(2-0.5,5-0.5) node {$\blue{2}$}; 
%\draw(2-0.5,6-0.5) node {$\blue{3}$};
%\draw(1-0.5,6-0.5) node {$\blue{4}$}; 
%\draw(1-0.5,7-0.5) node {$\blue{5}$};  
\end{tikzpicture}
}}\simeq \phi_{\ov52};\quad
%%%%
\theta_2=
\vcenter{\hbox{
\begin{tikzpicture}[x={(0in,-0.2in)},y={(0.2in,0in)}] % matrix coordinate
%shape phi
%\foreach \j in {6,...,7} \draw[thick] (1,\j) rectangle +(-1,-1);
%\foreach \j in {3,...,4} \draw[thick] (2,\j) rectangle +(-1,-1);
\foreach \j in {2,...,3} \draw[thick] (3,\j) rectangle +(-1,-1);
\foreach \j in {2,...,2} \draw[thick] (4,\j) rectangle +(-1,-1);
\foreach \j in {2,...,2} \draw[thick] (5,\j) rectangle +(-1,-1);
%content of phi
%\draw(5-0.5,1-0.5) node {$\green{\ov5}$};
\draw(5-0.5,2-0.5) node {$\green{\ov4}$};
\draw(4-0.5,2-0.5) node {$\green{\ov3}$};
\draw(3-0.5,2-0.5) node {$\green{\ov2}$}; 
\draw(3-0.5,3-0.5) node {$\green{\ov1}$};
%\draw(2-0.5,3-0.5) node {$\green{0}$}; 
%\draw(2-0.5,4-0.5) node {$\green{1}$}; 
%\draw(2-0.5,5-0.5) node {$\green{2}$}; 
%\draw(2-0.5,6-0.5) node {$\green{3}$};
%\draw(1-0.5,6-0.5) node {$\green{4}$}; 
%\draw(1-0.5,7-0.5) node {$\green{5}$};  
\end{tikzpicture}
}}\simeq \phi_{\ov4\ov1};\quad
%%%
\theta_3=
\vcenter{\hbox{
\begin{tikzpicture}[x={(0in,-0.2in)},y={(0.2in,0in)}] % matrix coordinate
%shape phi
%\foreach \j in {6,...,7} \draw[thick] (1,\j) rectangle +(-1,-1);
\foreach \j in {3,...,3} \draw[thick] (2,\j) rectangle +(-1,-1);
\foreach \j in {2,...,3} \draw[thick] (3,\j) rectangle +(-1,-1);
\foreach \j in {2,...,2} \draw[thick] (4,\j) rectangle +(-1,-1);
%\foreach \j in {2,...,2} \draw[thick] (5,\j) rectangle +(-1,-1);
%content of phi
%\draw(5-0.5,1-0.5) node {$\red{\ov5}$};
%\draw(5-0.5,2-0.5) node {$\red{\ov4}$};
\draw(4-0.5,2-0.5) node {$\red{\ov3}$};
\draw(3-0.5,2-0.5) node {$\red{\ov2}$}; 
\draw(3-0.5,3-0.5) node {$\red{\ov1}$};
\draw(2-0.5,3-0.5) node {$\red{0}$}; 
%\draw(2-0.5,4-0.5) node {$\red{1}$}; 
%\draw(2-0.5,5-0.5) node {$\red{2}$}; 
%\draw(2-0.5,6-0.5) node {$\red{3}$};
%\draw(1-0.5,6-0.5) node {$\red{4}$}; 
%\draw(1-0.5,7-0.5) node {$\red{5}$};  
\end{tikzpicture}
}}\simeq \phi_{\ov30};
\cr\cr
\theta_4=
\vcenter{\hbox{
\begin{tikzpicture}[x={(0in,-0.2in)},y={(0.2in,0in)}] % matrix coordinate
%shape phi
\foreach \j in {6,...,7} \draw[thick] (1,\j) rectangle +(-1,-1);
\foreach \j in {4,...,6} \draw[thick] (2,\j) rectangle +(-1,-1);
%\foreach \j in {2,...,3} \draw[thick] (3,\j) rectangle +(-1,-1);
%\foreach \j in {2,...,2} \draw[thick] (4,\j) rectangle +(-1,-1);
%\foreach \j in {1,...,2} \draw[thick] (5,\j) rectangle +(-1,-1);
%content of phi
%\draw(5-0.5,1-0.5) node {$\cyan{\ov5}$};
%\draw(5-0.5,2-0.5) node {$\cyan{\ov4}$};
%\draw(4-0.5,2-0.5) node {$\cyan{\ov3}$};
%\draw(3-0.5,2-0.5) node {$\cyan{\ov2}$}; 
%\draw(3-0.5,3-0.5) node {$\cyan{\ov1}$};
%\draw(2-0.5,3-0.5) node {$\cyan{0}$}; 
\draw(2-0.5,4-0.5) node {$\cyan{1}$}; 
\draw(2-0.5,5-0.5) node {$\cyan{2}$}; 
\draw(2-0.5,6-0.5) node {$\cyan{3}$};
\draw(1-0.5,6-0.5) node {$\cyan{4}$}; 
\draw(1-0.5,7-0.5) node {$\cyan{5}$};  
\end{tikzpicture}
}}\simeq \phi_{15};\quad
%%%%
\theta_5=
\vcenter{\hbox{
\begin{tikzpicture}[x={(0in,-0.2in)},y={(0.2in,0in)}] % matrix coordinate
%shape phi
%\foreach \j in {6,...,6} \draw[thick] (1,\j) rectangle +(-1,-1);
\foreach \j in {6,...,6} \draw[thick] (2,\j) rectangle +(-1,-1);
%\foreach \j in {2,...,3} \draw[thick] (3,\j) rectangle +(-1,-1);
%\foreach \j in {2,...,2} \draw[thick] (4,\j) rectangle +(-1,-1);
%\foreach \j in {1,...,2} \draw[thick] (5,\j) rectangle +(-1,-1);
%content of phi
%\draw(5-0.5,1-0.5) node {$\magenta{\ov5}$};
%\draw(5-0.5,2-0.5) node {$\magenta{\ov4}$};
%\draw(4-0.5,2-0.5) node {$\magenta{\ov3}$};
%\draw(3-0.5,2-0.5) node {$\magenta{\ov2}$}; 
%\draw(3-0.5,3-0.5) node {$\magenta{\ov1}$};
%\draw(2-0.5,3-0.5) node {$\magenta{0}$}; 
%\draw(2-0.5,4-0.5) node {$\magenta{1}$}; 
%\draw(2-0.5,5-0.5) node {$\magenta{2}$}; 
\draw(2-0.5,6-0.5) node {$\magenta{3}$};
%\draw(1-0.5,6-0.5) node {$\magenta{4}$}; 
%\draw(1-0.5,7-0.5) node {$\magenta{5}$};  
\end{tikzpicture}
}}\simeq \phi_{33};\qquad
\cr
\end{array}
\end{equation}

Two strips $\theta_p\simeq\phi_{ab}$ and $\theta_q\simeq\phi_{cd}$ from an outside decomposition 
may be combined to create a further strip via the $\#$ operation of Hamel and Goulden~\cite{HG}, 
which can be defined so that 
\begin{equation}\label{eqn-hash}
(\theta_p\#\theta_q)=\phi_{ad} \quad\mbox{where}\quad  \phi_{ad} \begin{cases} 
                        \mbox{is well defined for $a\leq d$};\cr
                        \mbox{is defined to be a null strip if $a=d+1$};\cr
                        \mbox{is empty if $a>d+1$}.\cr
                             \end{cases}
\end{equation}

Each diagram $(\theta_p\#\theta_q)\simeq\phi_{ad}$ is a substrip of the cutting strip $\phi$
extending from a box of content $a$ to one of content $d$. For $a\leq d$ its shape is that of the outer rim
of a skew diagram $F^{\nu/\kappa}$ for some partitions $\nu$ and $\kappa$. However,
the content of the boxes of $(\theta_p\#\theta_q)\simeq\phi_{ad}$ do not coincide with the conventionally 
defined content of the same boxes of $F^{\nu/\kappa}$, but instead they are shifted
by a fixed amount that we denote by $m(p,q)$ and refer to as the shift parameter of $(\theta_p\#\theta_q)$.
This is the value of the content of the box of $\phi_{ad}$ that lies at its intersection with the main 
diagonal of $F^\nu$.

In our running example the requisite hash products and shift parameters are given by:
\begin{equation}\label{eqn-mpq}
\begin{array}{|c|ccccc|}
\hline
(\theta_p\#\theta_q)&
&\theta_2&\theta_3&\theta_4&\theta_5\cr
\hline
&&&&&\cr
\theta_1&\phi_{\ov52}&\phi_{\ov5\ov1}&\phi_{\ov50}&\phi_{\ov55}&\phi_{\ov53}\cr
\theta_2&\phi_{\ov42}&\phi_{\ov4\ov1}&\phi_{\ov40}&\phi_{\ov45}&\phi_{\ov43}\cr
\theta_3&\phi_{\ov32}&\phi_{\ov3\ov1}&\phi_{\ov30}&\phi_{\ov35}&\phi_{\ov33}\cr
\theta_4&\phi_{12}&\phi_{1\ov1}&\phi_{10}&\phi_{15}&\phi_{13}\cr
\theta_5&\phi_{32}&\phi_{3\ov1}&\phi_{30}&\phi_{35}&\phi_{32}\cr
\hline
\end{array}
\qquad \qquad
\begin{array}{|c|ccccc|}
\hline
m(p,q)&\theta_1&\theta_2&\theta_3&\theta_4&\theta_5\cr
\hline
&&&&&\cr
\theta_1&\ov2&\ov3&\ov2&\ov1&\ov1\cr
\theta_2&\ov1&\ov2&\ov1&0&\ov1\cr
\theta_3&\ov1&\ov2&\ov1&0&\ov1\cr
\theta_4&1&-&0&2&1\cr
\theta_5&0&-&-&4&3\cr
\hline
\end{array}
\end{equation}

We are now in a position to state our main theorem.

\begin{Theorem}\label{The-HGXY} %[Hamel and Goulden] % outer decomposition identity}

For any partitions $\lambda$ and $\mu$ such that $\mu\subseteq\lambda$, 
and for all $\X=(x_{kc})$ and $\Y=(y_{\ell c})$, let the generalised skew 
Schur function $s_{\lambda/\mu}^{R'}(\X,\Y)$ be defined, 
as in Definition~\ref{Def-sfnXY-R'}, in terms of the skew supertableaux over the marked alphabet $R'$ 
obtained from $R=M\dot\cup N$ through Definition~\ref{Def-R'tableau}.
 
Let $\Theta=(\theta_1,\theta_2,\ldots,\theta_s)$ be the outside decomposition of $F^{\lambda/\mu}$
corresponding to a cutting strip $\phi$. Let $\tau$ be a shift operator whose action on  
$\X=(x_{kc})$ and $\Y=(y_{\ell c})$ is such that $\tau^t\X=(x_{k,t+c})$ and $\tau^t\Y=(y_{\ell,t+c})$
for any $t\in\Z$.

Then the generalised Schur functions of (\ref{eqn-sfnXY-R'}) satisfy the identity
\begin{equation}\label{eqn-HGX}
    s_{\lambda/\mu}^{R'}(\X,\Y) = \left|\,  s_{(\theta_p\#\theta_q)}^{R'}(\tau^{m(p,q)}\X,\tau^{m(p,q)}\Y) \,\right|_{1\leq p\leq q\leq s}  \,,
\end{equation}
where if $(\theta_p\#\theta_q)=\phi_{ad}$ then $m(p,q)$ is the value of the content of the box of $\phi_{ad}$ that 
lies at the intersection of the main diagonal of the diagram having $\phi_{ad}$ as its outer rim, and 
\begin{equation}\label{eqn-phashq}
   s_{(\theta_p\#\theta_q)}^{R'}(\tau^{m(p,q)}\X,\tau^{m(p,q)}\Y)= \begin{cases} 
	                                               s_{\phi_{ad}}^{R'}(\tau^{m(p,q)}\X,\tau^{m(p,q)}\Y)&\mbox{if $a\leq d$};\cr
		                                             1&\mbox{if $a=d+1$};\cr
																							   0&\mbox{if $a>d+1$}.\cr
															                \end{cases}
\end{equation}																	
\end{Theorem}

\noindent{\bf Proof}: We follow the proof of Theorem~\ref{The-HG} Hamel and Goulden~\cite{HG} which used the technique of Gessel-Viennot-Lindstr\"om  \cite{GV}, \cite{Lin} to construct a bijection between semistandard skew tableaux and non-intersecting lattice paths. In their construction Hamel and Goulden situate the paths on a rectangular lattice whose directed edges are defined by the orientation of the strips in the outside decomposition. Using cutting strip language this means the lattice has diagonal or horizontal edges in certain positions depending on whether a box in the cutting strip has a box below or beside it.  Details are made precise below in the context of our supertableaux.

We construct the lattice as a coordinate grid where the directed edges are determined
by 
our choice of alphabet $R'$ and cutting strip $\phi$ which themselves determine the sets $M$, $N$, $L$ and $B$. 
The lattice is constructed by noting that each edge terminating at coordinates $(r,c)$ is specified as follows:
\begin{equation}\label{eqn-edges}
\begin{array}{|c|c|c|c|c|c|}
\hline
\blue{r\in M},\brown{c\in L}&\red{r\in N},\green{c\in B}&\blue{r\in M},\green{c\in B}&\red{r\in N},\brown{c\in L}&\brown{c+1\in L}&\green{c+1\in B}\cr
\hline
&&&&&\cr
\vcenter{\hbox{
\begin{tikzpicture}[x={(0in,-0.4in)},y={(0.4in,0in)}] 
%horizontal left to right edges
\draw[-latex,blue,thick](1,2-1)--(1,2-0.4); \draw[blue,thick](1,2-0.4)--(1,2); \draw(1,2)node{$\bullet$};
\end{tikzpicture}
}}
&
\vcenter{\hbox{
\begin{tikzpicture}[x={(0in,-0.4in)},y={(0.4in,0in)}] 
%horizontal left to right edges
\draw[-latex,red,thick](1,2-1)--(1,2-0.4); \draw[red,thick](1,2-0.4)--(1,2); \draw(1,2)node{$\bullet$};
\end{tikzpicture}
}}
&
%diagonal up edge
\vcenter{\hbox{
\begin{tikzpicture}[x={(0in,-0.4in)},y={(0.4in,0in)}] 
\draw[-latex,blue,thick](1+1,2-1)--(1+0.4,2-0.4); \draw[blue,thick](1+0.4,2-0.4)--(1,2); \draw(1,2)node{$\bullet$};
\end{tikzpicture}
}}
&
%diagonal down edges
\vcenter{\hbox{
\begin{tikzpicture}[x={(0in,-0.4in)},y={(0.4in,0in)}] 
\draw[-latex,red,thick](1-1,2-1)--(1-0.4,2-0.4); \draw[red,thick](1-0.4,2-0.4)--(1,2); \draw(1,2)node{$\bullet$};
\end{tikzpicture}
}}
&
%vertical down edges
\vcenter{\hbox{
\begin{tikzpicture}[x={(0in,-0.4in)},y={(0.4in,0in)}] 
\draw[-latex,thick](1-1,2)--(1-0.4,2); \draw[thick](1-0.4,2)--(1,2); \draw(1,2)node{$\bullet$};
\end{tikzpicture}
}}
&
%vertical up edges
\vcenter{\hbox{
\begin{tikzpicture}[x={(0in,-0.4in)},y={(0.4in,0in)}] 
\draw[-latex,thick](1+1,2)--(1+0.4,2); \draw[thick](1+0.4,2)--(1,2); \draw(1,2)node{$\bullet$};
\end{tikzpicture}
}}
\cr
&&&&&\cr
\hline
\end{array}
\end{equation}

In our example it follows that the required lattice takes the form
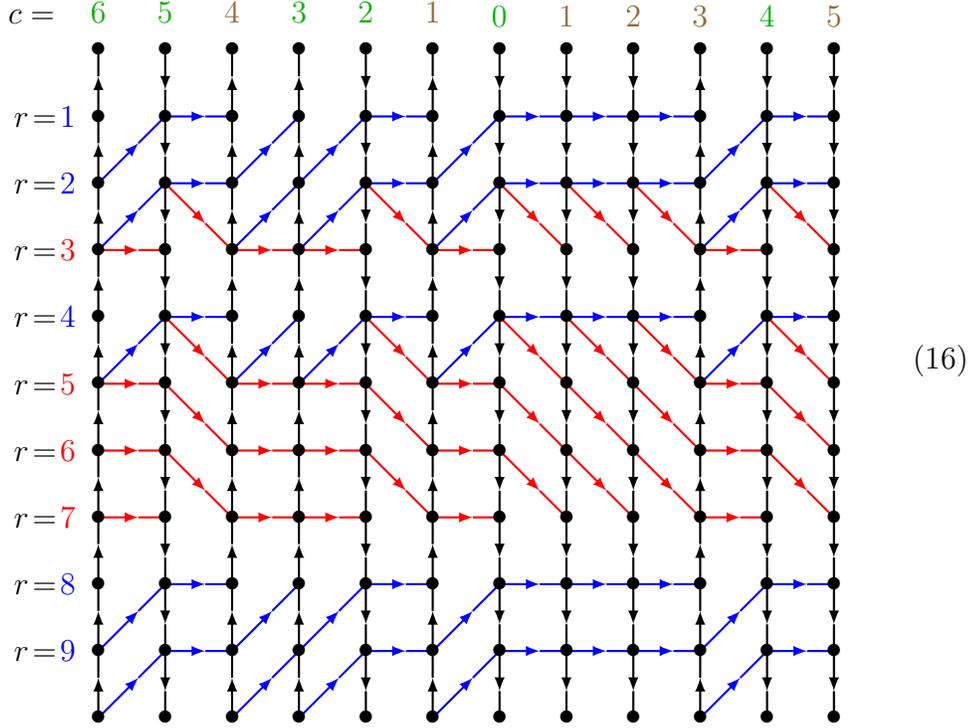
\begin{equation}\label{eqn-lattice}
\vcenter{\hbox{
\begin{tikzpicture}[x={(0in,-0.35in)},y={(0.35in,0in)}] %[x={(0in,-0.4in)},y={(0.4in,0in)}] % matrix coordinates
%coordinate labels
%\draw(0.5,-5)node{$k$};
\foreach \i in {1,2,4,8,9} \draw(\i,-6.8)node{$r\!=\!\blue{\i}$};
\foreach \i in {3,5,6,7} \draw(\i,-6.8)node{$r\!=\!\red{\i}$};
\draw(-0.5,-7)node{$c=$};
\foreach \j in {2,3,5,6} \draw(-0.6,-\j)node{$\green{\ov\j}$};
\foreach \j in {1,4} \draw(-0.6,-\j)node{$\brown{\ov\j}$};
\foreach \j in {0,4} \draw(-0.5,\j)node{$\green{\j}$};
\foreach \j in {1,2,3,5} \draw(-0.5,\j)node{$\brown{\j}$};
%horizontal left to right edges
\foreach \i in {1,2,4,8,9} \foreach \j in {-4,-1,1,2,3,5} \draw[-latex,blue,thick](\i,\j-1)--(\i,\j-0.4);
\foreach \i in {1,2,4,8,9} \foreach \j in {-4,-1,1,2,3,5} \draw[blue,thick](\i,\j-0.4)--(\i,\j);
\foreach \i in {3,5,6,7} \foreach \j in {-5,-3,-2,0,4} \draw[-latex,red,thick](\i,\j-1)--(\i,\j-0.4);
\foreach \i in {3,5,6,7} \foreach \j in {-5,-3,-2,0,4} \draw[red,thick](\i,\j-0.4)--(\i,\j);
%diagonal up edge
\foreach \i in {1,2,4,8,9} \foreach \j in {-5,-3,-2,0,4} \draw[-latex,blue,thick](\i+1,\j-1)--(\i+0.4,\j-0.4); %diagonally up arrows [-latex]
\foreach \i in {1,2,4,8,9} \foreach \j in {-5,-3,-2,0,4} \draw[blue,thick](\i+0.4,\j-0.4)--(\i,\j);
%diagonal down edges
\foreach \i in {3,5,6,7} \foreach \j in {-4,-1,1,2,3,5} \draw[-latex,red,thick](\i-1,\j-1)--(\i-0.4,\j-0.4); %diagonally down arrows [-latex]
\foreach \i in {3,5,6,7} \foreach \j in {-4,-1,1,2,3,5} \draw[red,thick](\i-0.4,\j-0.4)--(\i,\j);
%vertical down edges
\foreach \j in {-5,-2,0,1,2,4,5} \foreach \i in {1,...,10} \draw[-latex,thick](\i-1,\j)--(\i-0.4,\j);
\foreach \j in {-5,-2,0,1,2,4,5} \foreach \i in {1,...,10}\draw[thick](\i-0.4,\j)--(\i,\j);
%vertical up edges
\foreach \j in {-6,-4,-3,-1,3} \foreach \i in {0,...,9} \draw[-latex,thick](\i+1,\j)--(\i+0.4,\j);
\foreach \j in {-6,-4,-3,-1,3} \foreach \i in {0,...,9} \draw[thick](\i+0.4,\j)--(\i,\j);
%path
%\draw[draw=blue,ultra thick] (0,-3)--(3,-3)--(3,-2)--(2,-2)--(1,-1)--(1,0)--(2,0)--(2,1)--(3,1)--(3,2)--(2,3)--(4,3)--(4,4)--(5,4);
%lattice grid
\foreach \i in {0,...,10} \foreach \j in {-6,...,5} \draw(\i,\j)node{$\bullet$};
%label P and Q nodes
%\draw(5,-4.4)node{$\blue{P_1}$}; \draw(5,-4)node{$\blue\bullet$}; \draw(5,0.6)node{$\blue{Q_1}$}; \draw(5,1)node{$\blue\bullet$};
%\draw(0,-3.4)node{$\red{P_2}$}; \draw(0,-3)node{$\red\bullet$}; \draw(0,-2.4)node{$\red{Q_2}$}; \draw(0,-2)node{$\red\bullet$};
%\draw(0,-0.4)node{$\cyan{P_3}$}; \draw(0,0)node{$\cyan\bullet$}; \draw(5,3.6)node{$\cyan{Q_3}$}; \draw(5,4)node{$\cyan\bullet$};
\end{tikzpicture}
}}
\end{equation}

As in \cite{HG} for each strip $\theta_p$ in $T$  we can determine the start points $P_p$ and the end points $Q_p$ of the corresponding path. These are exactly determined by the location of $\phi_{ab}\simeq\theta_p$ in $F^{\lambda/\mu}$:
\begin{itemize}
\item if $\theta_p\simeq\phi_{ab}$ starts on the  left boundary then $P_p=(0,a-1)$;
\item if $\theta_p\simeq\phi_{ab}$ starts on the bottom boundary then $P_p=(n+1,a-1)$;
\item if $\theta_p\simeq\phi_{ab}$ ends on  the right boundary then  $Q_p=(n+1,b)$;
\item if $\theta_p\simeq\phi_{ab}$ ends on the top boundary then $Q_p=(0,b)$.
\end{itemize}

The bijection is defined from a tableau $T$ of shape $F^{\lambda/\mu}$ to an 
$s$-tuple of lattice paths with start points and  end points determined by
the outside decomposition $\Theta$. This map  sends an 
entry $r$ or $r'$ in a box of content $c$ in $T$
to a horizontal or diagonal edge whose right hand endpoint is situated at $(r,c)$ in accordance with the
prescription (\ref{eqn-edges}).  All that remains is to specify the vertical edges. These are simply the edges required to complete the path in a continuous manner, and these vertical edges are again as prescribed by (\ref{eqn-edges}).

The process can be exemplified as follows. 
In the case $\lambda=(6,5,5,4,4,3)$ and $\mu=(3,2)$ a typical supertableau $T\in{\cal T}_{R'}^{\lambda/\mu}$ with entries $k$ from 
$\{1<2<3'<4<5'<6'<7'<8<9\}$ is as shown below on the right alongside a reminder of the content $c$ of the boxes of $F^{\lambda/\mu}$ 
on the left.
\begin{equation}\label{eqn-FT}
\begin{array}{ccc}
\mbox{Content of $F^{\lambda/\mu}$}&\qquad\qquad&\mbox{Tableau $T^{\lambda/\mu}$}\cr
%\cr
\vcenter{\hbox{
\begin{tikzpicture}[x={(0in,-0.2in)},y={(0.2in,0in)}]%[x={(0in,-0.2in)},y={(0.2in,0in)}] % matrix coordinate
%shape lambda
\foreach \j in {1,...,6} \draw[thick] (1,\j) rectangle +(-1,-1); 
\foreach \j in {1,...,5} \draw[thick] (2,\j) rectangle +(-1,-1);
\foreach \j in {1,...,5} \draw[thick] (3,\j) rectangle +(-1,-1);
\foreach \j in {1,...,4} \draw[thick] (4,\j) rectangle +(-1,-1);
\foreach \j in {1,...,4} \draw[thick] (5,\j) rectangle +(-1,-1);
\foreach \j in {1,...,3} \draw[thick] (6,\j) rectangle +(-1,-1);
%shape mu
\foreach \j in {1,...,3} \draw(1-0.5,\j-0.5)node{$\ast$};
\foreach \j in {1,...,2} \draw(2-0.5,\j-0.5)node{$\ast$};
%content of lambda/mu
\foreach \k in {3,...,5} \draw(1-0.5,\k+0.5) node{$\k$};
\foreach \k in {1,...,3} \draw(2-0.5,\k+1.5) node{$\k$};
\foreach \k in {1,...,2} \draw(3-0.5,-\k+2.5) node{$\ov{\k}$};\foreach \k in {0,...,2} \draw(3-0.5,\k+2.5) node{$\k$};
\foreach \k in {1,...,3} \draw(4-0.5,-\k+3.5) node{$\ov{\k}$};\foreach \k in {0,...,0} \draw(4-0.5,\k+3.5) node{$\k$};
\foreach \k in {1,...,4} \draw(5-0.5,-\k+4.5) node{$\ov{\k}$};
\foreach \k in {3,...,5} \draw(6-0.5,-\k+5.5) node{$\ov{\k}$};
\end{tikzpicture}
}}
&&
\vcenter{\hbox{
\begin{tikzpicture}[x={(0in,-0.2in)},y={(0.2in,0in)}] % matrix coordinate
%shape lambda
\foreach \j in {1,...,6} \draw[thick] (1,\j) rectangle +(-1,-1); 
\foreach \j in {1,...,5} \draw[thick] (2,\j) rectangle +(-1,-1);
\foreach \j in {1,...,5} \draw[thick] (3,\j) rectangle +(-1,-1);
\foreach \j in {1,...,4} \draw[thick] (4,\j) rectangle +(-1,-1);
\foreach \j in {1,...,4} \draw[thick] (5,\j) rectangle +(-1,-1);
\foreach \j in {1,...,3} \draw[thick] (6,\j) rectangle +(-1,-1);
%shape mu
\foreach \j in {1,...,3} \draw(1-0.5,\j-0.5)node{$\ast$};
\foreach \j in {1,...,2} \draw(2-0.5,\j-0.5)node{$\ast$};
%content of lambda/mu
\draw(1-0.5,3+0.5) node{$\magenta{3'}$};\draw(1-0.5,4+0.5) node{$\cyan{5'}$};\draw(1-0.5,5+0.5) node{$\cyan{8}$};
\draw(2-0.5,2+0.5) node{$\cyan{2}$};\draw(2-0.5,3+0.5) node{$\cyan{3'}$};\draw(2-0.5,4+0.5) node{$\cyan{5'}$};
\draw(3-0.5,0+0.5) node{$\green{1}$};\draw(3-0.5,1+0.5) node{$\green{1}$};\draw(3-0.5,2+0.5) node{$\blue{3'}$};\draw(3-0.5,3+0.5) node{$\blue{4}$};\draw(3-0.5,4+0.5) node{$\blue{9}$};
\draw(4-0.5,0+0.5) node{$\green{3'}$};\draw(4-0.5,1+0.5) node{$\blue{4}$};\draw(4-0.5,2+0.5) node{$\blue{4}$};\draw(4-0.5,3+0.5) node{$\red{6'}$};
\draw(5-0.5,0+0.5) node{$\green{3'}$};\draw(5-0.5,1+0.5) node{$\blue{7'}$};\draw(5-0.5,2+0.5) node{$\red{8}$};\draw(5-0.5,3+0.5) node{$\red{8}$};
\draw(6-0.5,0+0.5) node{$\blue{4}$};\draw(6-0.5,1+0.5) node{$\blue{8}$};\draw(6-0.5,2+0.5) node{$\red{9}$};
\end{tikzpicture}
}}\cr
\end{array}
\end{equation}

The corresponding $5$-tuple of non-intersecting lattice paths is illustrated below, where for the sake of clarity the directions 
of edges have been omitted:
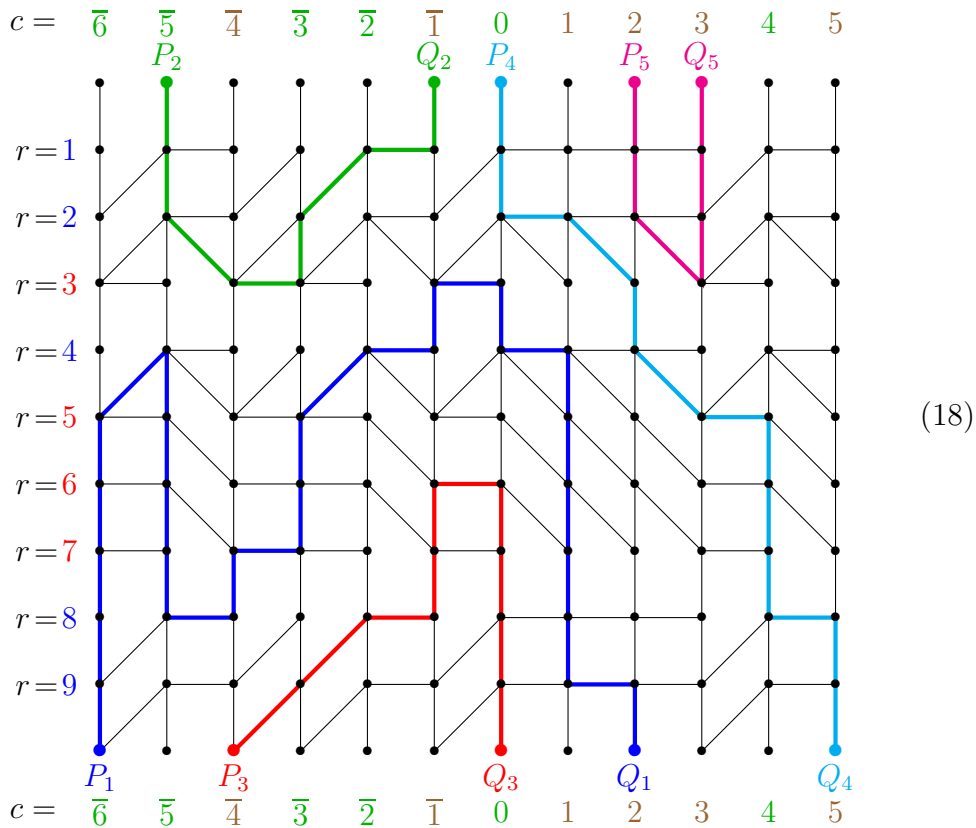
\begin{equation}
\vcenter{\hbox{
\begin{tikzpicture}[x={(0in,-0.35in)},y={(0.35in,0in)}] %[x={(0in,-0.4in)},y={(0.4in,0in)}] % matrix coordinates
%horizontal left to right edges
\foreach \i in {1,2,4,8,9} \foreach \j in {-4,-1,1,2,3,5} \draw(\i,\j-1)--(\i,\j);
\foreach \i in {3,5,6,7} \foreach \j in {-5,-3,-2,0,4} \draw(\i,\j-1)--(\i,\j);
%diagonal up edge
\foreach \i in {1,2,4,8,9} \foreach \j in {-5,-3,-2,0,4} \draw(\i+1,\j-1)--(\i,\j);
%diagonal down edges
\foreach \i in {3,5,6,7} \foreach \j in {-4,-1,1,2,3,5} \draw(\i-1,\j-1)--(\i,\j);
%vertical down edges
\foreach \j in {-5,-2,0,1,2,4,5} \foreach \i in {1,...,10}\draw(\i-1,\j)--(\i,\j);
%vertical up edges
\foreach \j in {-6,-4,-3,-1,3} \foreach \i in {0,...,9} \draw(\i+1,\j)--(\i,\j);
%path
\draw[draw=blue,ultra thick] (10,-6)--(5,-6)--(4,-5)--(8,-5)--(8,-4)--(7,-4)--(7,-3)--(5,-3)--(4,-2)--(4,-1)--(3,-1)--(3,0)--(4,0)--(4,1)--(9,1)--(9,2)--(10,2);
\draw[draw=green!70!black,ultra thick] (0,-5)--(2,-5)--(3,-4)--(3,-3)--(2,-3)--(1,-2)--(1,-1)--(0,-1);
\draw[draw=red,ultra thick] (10,-4)--(8,-2)--(8,-1)--(6,-1)--(6,-0)--(10,0);
\draw[draw=cyan,ultra thick] (0,0)--(2,0)--(2,1)--(3,2)--(4,2)--(5,3)--(5,4)--(8,4)--(8,5)--(10,5);
\draw[draw=magenta,ultra thick] (0,2)--(2,2)--(3,3)--(0,3);
%lattice grid
\foreach \i in {0,...,10} \foreach \j in {-6,...,5} \draw(\i,\j)node{$\sc\bullet$};
%label P and Q nodes
\draw(10.4,-6)node{$\blue{P_1}$}; \draw(10,-6)node{$\blue\bullet$}; \draw(10.4,2)node{$\blue{Q_1}$}; \draw(10,2)node{$\blue\bullet$};
\draw(-0.4,-5)node{$\green{P_2}$}; \draw(0,-5)node{$\green\bullet$}; \draw(-0.4,-1)node{$\green{Q_2}$}; \draw(0,-1)node{$\green\bullet$};
\draw(10.4,-4)node{$\red{P_3}$}; \draw(10,-4)node{$\red\bullet$}; \draw(10.4,0)node{$\red{Q_3}$}; \draw(10,0)node{$\red\bullet$};
\draw(-0.4,0)node{$\cyan{P_4}$}; \draw(0,0)node{$\cyan\bullet$}; \draw(10.4,5)node{$\cyan{Q_4}$}; \draw(10,5)node{$\cyan\bullet$};
\draw(-0.4,2)node{$\magenta{P_5}$}; \draw(0,2)node{$\magenta\bullet$}; \draw(-0.4,3)node{$\magenta{Q_5}$}; \draw(0,3)node{$\magenta\bullet$};
%coordinate labels (k,c)
%node{$k=$};
\foreach \i in {1,2,4,8,9} \draw(\i,-6.8)node{$r\!=\!\blue{\i}$};
\foreach \i in {3,5,6,7} \draw(\i,-6.8)node{$r\!=\!\red{\i}$};
%node{$c=$}
\draw(-0.9,-7)node{$c=$};
\foreach \j in {2,3,5,6} \draw(-0.9,-\j)node{$\green{\ov\j}$};
\foreach \j in {1,4} \draw(-0.9,-\j)node{$\brown{\ov\j}$};
\foreach \j in {0,4} \draw(-0.9,\j)node{$\green{\j}$};
\foreach \j in {1,2,3,5} \draw(-0.9,\j)node{$\brown{\j}$};
\draw(10.9,-7)node{$c=$};
\foreach \j in {2,3,5,6} \draw(10.9,-\j)node{$\green{\ov\j}$};
\foreach \j in {1,4} \draw(10.9,-\j)node{$\brown{\ov\j}$};
\foreach \j in {0,4} \draw(10.9,\j)node{$\green{\j}$};
\foreach \j in {1,2,3,5} \draw(10.9,\j)node{$\brown{\j}$};
\end{tikzpicture}
}}
\end{equation}

Hamel \cite{Ham} provided a bijection between shifted tableaux for Schur $Q$-functions and tuples 
of nonintersecting lattice paths.  Like our supertableaux introduced 
in Definition \ref{Def-R'tableau}, the shifted tableaux have entries that weakly increase in rows and columns  and  
have at most one primed entry per row and at most one unprimed entry per column. Our outer boundary is different from that of the shifted tableaux and the order of the entries may also be different since the order adopted in~\cite{Ham} corresponds to setting
$R'=\{1'<1<2'<2< \cdots\}$.  The outer boundary accounts for the placement of the starting and ending points, $P$s and $Q$s; however,  the argument in 
\cite{Ham} is not specific to that context and would cover the placement of $P$s and $Q$s defined here. Moreover, the order of entries is not an issue as the entries from $R'$ are still organized in the canonical counting order of $R$ even if some can only appear in either primed or unprimed incarnations.
This difference does not affect the argument that the image of supertableaux under our bijective map remain
tuples of nonintersecting lattice paths.

Now we may use the weight preserving map of
\cite{Ham} to associate
each entry  $k$ or $k'$ in the supertableau to a horizontal of diagonal directed edge with 
right hand endpoint at level
corresponding to $k$ or $k'$ in the $s$-tuple of non-intersecting
lattice paths.  
In the simplest first variation context
both the supertableau entry and the lattice path edge carry
a weight of $x_k$.  The vertical edges in the lattice path are considered to have weight $1$.
As in \cite{FK20}, this identification of entry and edge can be weighted more generally
to give the ninth variation. 
In this case, with a fixed but arbitrary choice of $\X=(x_{kc})$ and $\Y=(y_{\ell c})$, we assign 
weights $x_{kc}$ and $y_{\ell c}$ to each entry $r=i_k\in M$ or $r'=j'_\ell\in N'$ on the $c$th diagonal of the supertableau,
and we use the same bijection as before, mapping the
weighted supertableau entry to a horizontal or diagonal directed edge
in the vertical swath between vertical lines labelled by $c-1$ and $c$, assigning an identical weight to this edge.

Then Stembridge's interpretation~\cite{Ste90} of Gessel-Viennot-Lindstr\"om~\cite{GV} \cite{Lin} with weights 
and underlying lattice (called {\em digraph} in \cite{Ste90}) as described above and exemplified in 
(\ref{eqn-lattice}) completes the proof of Theorem~\ref{The-HGXY}.
\qed

\section{Factorial supersymmetric skew Schur functions}\label{Sec:factssfn}

The specialisation we consider here is one in which $s_{\lambda/\mu}^{R'}(\X,\Y)$ is of the form
$s^{R'}_{\lambda/\mu}(\x,\y|\a)$ where $\x=(x_1,x_2,\ldots,x_m)$, $\y=(y_1,y_2,\ldots,n)$
and $\a$ is a sequence of factorial parameters. The aim is to provide 
a definition such that $s^{R'}_{\lambda/\mu}(\x,\y|\a)$ is supersymmetric	for all $\a$, that is to say is
symmetric under independent arbitrary permutations of the components of $\x$ and $\y$ and is independent of
$t$ if $x_k=t=-y_\ell$ for any $k$ and $\ell$. As a fortunate by-product of this definition it will be shown 
that $s^{R'}_{\lambda/\mu}(\x,\y|\a)$ is also independent of $R'=M\dot\cup N'$ in the sense that it does not depend on 
the choice of the elements of $M$ and $N$, only on their sizes $m$ and $n$, respectively. 
Prior to making our main definition it is helpful to introduce some further notation relating to the 
specification of $R'$ as follows:

\begin{Definition}\label{Def-sigma}
For $R=M\dot\cup N$ with $M=\{i_1<i_2<\cdots<i_m\}$ and $N=\{j_1<j_2<\cdots<j_n\}$ let
\begin{equation}
    \sigma(r)=\begin{cases}
		                \#\{i\in M\,|\,i\leq r\}-\#\{j\in N\,|\,j\leq r\}&\mbox{if $r\in M$};\cr
										\#\{i\in M\,|\,i\leq r\}-\#\{j\in N\,|\,j\leq r\}+1&\mbox{if $r\in N$}.\cr
		            \end{cases}
\end{equation}
\end{Definition}

For example, if $R=\{1<2<\cdots<9\}$ with $M=\{1<2<4<8<9\}$ and $N=\{3<5<6<7\}$ then we have
\begin{equation}
\begin{array}{|c||c|c|c|c|c|c|c|c|c|}
\hline
r=&1&2&3&4&5&6&7&8&9\cr
\hline
\sigma(r)=&1&2&2&2&2&1&0&0&1\cr
\hline
\end{array}
\end{equation}

With this notation our factorial supersymmetric skew Schur functions are given by:
\begin{Definition}\label{Def-fact-sfnxy}
For partitions $\lambda$ and $\mu$ such that $\nu\subseteq\lambda$ let
${\cal T}^{\lambda/\mu}_{R'}$ be the set of all skew supertableaux of shape $F^{\lambda/\mu}$
with entries from the alphabet $R'$ obtained from $R=(1<2<\cdots<m+n)$ by
leaving elements of the subalphabet $M=(i_1<i_2<\cdots<i_m)$ unprimed while adding primes 
to elements of the complementary subalphabet $N=(j_1<j_2<\cdots<j_n)$.
Then for all $\x=(x_1,x_2,\ldots,x_m)$, $\y=(y_1,y_2,\ldots,y_n)$ and $\a=(\ldots,a_{\ov2},a_{\ov1},a_0,a_1,a_2,\ldots)$ 
let
\begin{equation}
   s^{R'}_{\lambda/\mu}(\x,\y|\a)= \sum_{T\in{\cal T}^{\lambda/\mu}_{R'}} \prod_{(i,j)\in F^{\lambda/\mu}}  \wgt(t_{ij})
\end{equation} 
where $t_{ij}$ is the entry in the box in $i$th row and $j$th column of $T$, with 
\begin{equation}
               \wgt(t_{ij}) = \begin{cases}
							                     x_k+a_{\sigma(r)+j-i}&\mbox{if $t_{ij}=r=i_k\in M$};\cr
																	 y_\ell-a_{\sigma(r)+j-i}&\mbox{if $t_{ij}=r'=j'_\ell\in N'$}.\cr
							\end{cases}
\end{equation}
\end{Definition}

For $T$ as in (\ref{eqn-tab-ex}) with $R'=\{\blue1<\blue2<\red{3'}<\blue4<\red{5'}<\red{6'}<\red{7'}<\blue8<\blue9\}$
the corresponding contribution to $s^{R'}_{\lambda/\mu}(\x,\y|\a)$ 
is the product of the terms in the following array:
\begin{equation}
\vcenter{\hbox{
\begin{tikzpicture}[x={(0in,-0.2in)},y={(0.6in,0in)}] % matrix coordinate
%shape lambda
\foreach \j in {1,...,6} \draw[thick] (1,\j) rectangle +(-1,-1); 
\foreach \j in {1,...,5} \draw[thick] (2,\j) rectangle +(-1,-1);
\foreach \j in {1,...,5} \draw[thick] (3,\j) rectangle +(-1,-1);
\foreach \j in {1,...,4} \draw[thick] (4,\j) rectangle +(-1,-1);
\foreach \j in {1,...,4} \draw[thick] (5,\j) rectangle +(-1,-1);
\foreach \j in {1,...,3} \draw[thick] (6,\j) rectangle +(-1,-1);
%shape mu
\foreach \j in {1,...,3} \draw(1-0.5,\j-0.5)node{$\ast$};
\foreach \j in {1,...,2} \draw(2-0.5,\j-0.5)node{$\ast$};
%content of lambda/mu
%content of lambda/mu
\draw(1-0.5,3+0.5) node{$\red{y_{1}-a_{5}}$};\draw(1-0.5,4+0.5) node{$\red{y_{2}-a_{6}}$};\draw(1-0.5,5+0.5) node{$\blue{x_{4}+a_{5}}$};
\draw(2-0.5,2+0.5) node{$\blue{x_{2}+a_{3}}$};\draw(2-0.5,3+0.5) node{$\red{y_{1}-a_{4}}$};\draw(2-0.5,4+0.5) node{$\red{y_{2}-a_{5}}$};
\draw(3-0.5,0+0.5) node{$\blue{x_{1}+a_{\ov1}}$};\draw(3-0.5,1+0.5) node{$\blue{x_{1}+a_{0}}$};\draw(3-0.5,2+0.5) node{$\red{y_{1}-a_{2}}$};\draw(3-0.5,3+0.5) node{$\blue{x_{3}+a_{3}}$};\draw(3-0.5,4+0.5) node{$\blue{x_{5}+a_{3}}$};
\draw(4-0.5,0+0.5) node{$\red{y_{1}-a_{\ov1}}$};\draw(4-0.5,1+0.5) node{$\blue{x_{3}+a_{0}}$};\draw(4-0.5,2+0.5) node{$\blue{x_{3}+a_{1}}$};\draw(4-0.5,3+0.5) node{$\red{y_{3}-a_{1}}$};
\draw(5-0.5,0+0.5) node{$\red{y_{1}-a_{\ov2}}$};\draw(5-0.5,1+0.5) node{$\red{y_{4}-a_{\ov3}}$};\draw(5-0.5,2+0.5) node{$\blue{x_{4}+a_{\ov2}}$};\draw(5-0.5,3+0.5) node{$\blue{x_{4}+a_{\ov1}}$};
\draw(6-0.5,0+0.5) node{$\blue{x_{3}+a_{\ov3}}$};\draw(6-0.5,1+0.5) node{$\blue{x_{4}+a_{\ov4}}$};\draw(6-0.5,2+0.5) node{$\blue{x_{5}+a_{\ov2}}$};
\end{tikzpicture}
}} \label{tableau}
\end{equation}

The map of each entry of $T$ to a node within our lattice, along with its corresponding weight, takes the
following form in which rather than displaying any lattice paths or even lattice edges these
are omitted in favour of lines of constant $r$ and constant $\sigma(r)+c$. 
For $r\in M$ or $r\in N$ the lines of constant $r$ are horizontal $x$-lines or $y$-lines labelled by 
parameters $x_k$ and $y_\ell$, respectively. The lines of constant $\sigma(r)+c$, 
labelled by $a_{\sigma(r)+c}$, slope to the left, or are vertical or slope to the right
according as they pass downwards between a pair of $x$-lines,  or 
an $x$-line and a $y$-line or {\it vice versa}, or a pair of $y$-lines, respectively.
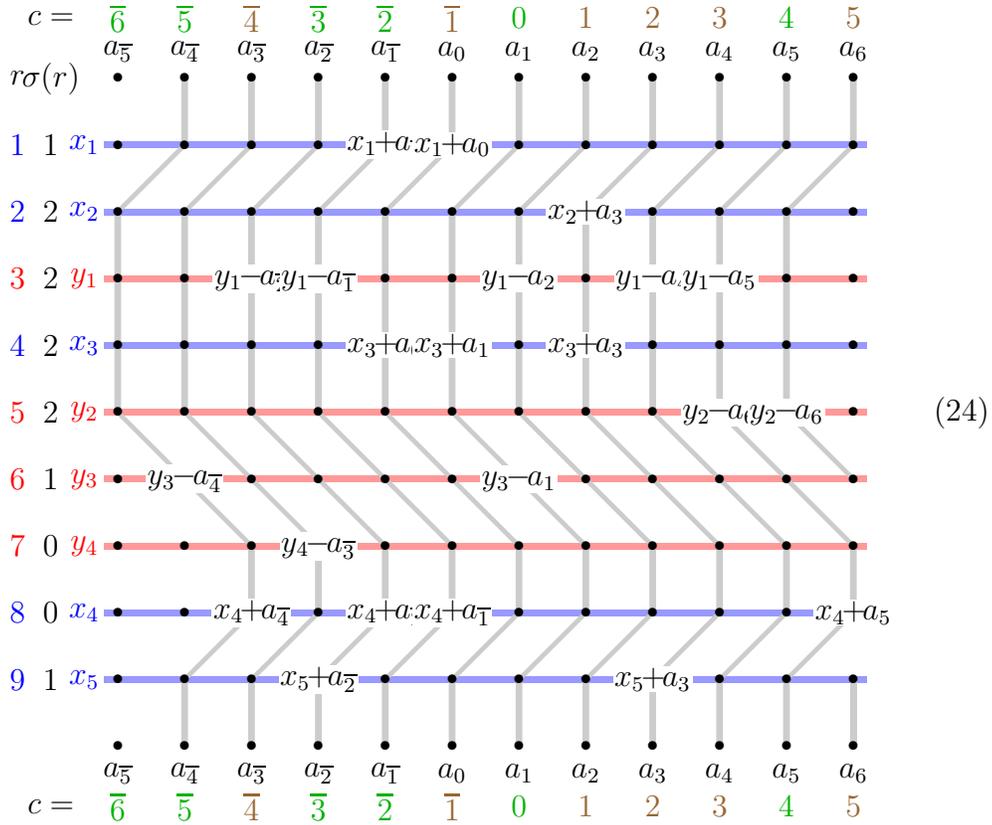
\begin{equation}
\vcenter{\hbox{
\begin{tikzpicture}[x={(0in,-0.35in)},y={(0.35in,0in)}]%[x={(0in,-0.4in)},y={(0.45in,0in)}] % matrix coordinates
%horizontal x_k and y_l 
\foreach \i in {1,2,4,8,9} \foreach \j in {0.05} \fill[blue!40!white](\i-\j,-6.2)--(\i-\j,5.2)--(\i+\j,5.2)--(\i+\j,-6.2)--cycle;
\foreach \i in {3,5,6,7} \foreach \j in {0.05} \fill[red!40!white](\i-\j,-6.2)--(\i-\j,5.2)--(\i+\j,5.2)--(\i+\j,-6.2)--cycle;
%paths of constant a_m
\foreach \i in {-4,...,5} \foreach \j in {0.05} \fill[black!20!white]
(0,\i-1+\j)--(1,\i-1+\j)--(2,\i-2+\j)--(3,\i-2+\j)--(4,\i-2+\j)--(5,\i-2+\j)--(6,\i-1+\j)--(7,\i-0+\j)--(8,\i-0+\j)--(9,\i-1+\j)--(10,\i-1+\j)
--(10,\i-1-\j)--(9,\i-1-\j)--(8,\i-0-\j)--(7,\i-0-\j)--(6,\i-1-\j)--(5,\i-2-\j)--(4,\i-2-\j)--(3,\i-2-\j)--(2,\i-2-\j)--(1,\i-1-\j)--(0,\i-1-\j)--cycle;
\foreach \i in {6} \foreach \j in {0.05} \fill[black!20!white]
(0,\i-1+\j)--(1,\i-1+\j)--(2,\i-2+\j)--(3,\i-2+\j)--(4,\i-2+\j)--(5,\i-2+\j)--(6,\i-1+\j)
%--(7,\i-0+\j)--(8,\i-0+\j)--(9,\i-1+\j)--(10,\i-1+\j)
%--(10,\i-1-\j)--(9,\i-1-\j)--(8,\i-0-\j)--(7,\i-0-\j)
--(6,\i-1-\j)--(5,\i-2-\j)--(4,\i-2-\j)--(3,\i-2-\j)--(2,\i-2-\j)--(1,\i-1-\j)--(0,\i-1-\j)--cycle;
\foreach \i in {6} \foreach \j in {0.05} \fill[black!20!white](9,\i-1+\j)--(10,\i-1+\j)--(10,\i-1-\j)--(9,\i-1-\j)--cycle;
%path
%\draw[draw=blue,ultra thick] (10,-6)--(5,-6)--(4,-5)--(8,-5)--(8,-4)--(7,-4)--(7,-3)--(5,-3)--(4,-2)--(4,-1)--(3,-1)--(3,0)--(4,0)--(4,1)--(9,1)--(9,2)--(10,2);
%\draw[draw=green!70!black,ultra thick] (0,-5)--(2,-5)--(3,-4)--(3,-3)--(2,-3)--(1,-2)--(1,-1)--(0,-1);
%\draw[draw=red,ultra thick] (10,-4)--(8,-2)--(8,-1)--(6,-1)--(6,-0)--(10,0);
%\draw[draw=cyan,ultra thick] (0,0)--(2,0)--(2,1)--(3,2)--(4,2)--(5,3)--(5,4)--(8,4)--(8,5)--(10,5);
%\draw[draw=magenta,ultra thick] (0,2)--(2,2)--(3,3)--(0,3);
%lattice grid
\foreach \i in {0,...,10} \foreach \j in {-6,...,5} \draw(\i,\j)node{$\sc\bullet$};
%weights
\draw(1,-2)node[fill=white,inner sep=0.5pt]{$x_1\!\!+\!\!a_{\ov1}$};
\draw(1,-1)node[fill=white,inner sep=0.5pt]{$x_1\!\!+\!\!a_{0}$};
\draw(2,1)node[fill=white,inner sep=0.5pt]{$x_2\!\!+\!\!a_{3}$};
\draw(3,-4)node[fill=white,inner sep=0.5pt]{$y_1\!\!-\!\!a_{\ov2}$};
\draw(3,-3)node[fill=white,inner sep=0.5pt]{$y_1\!\!-\!\!a_{\ov1}$};
\draw(3,0)node[fill=white,inner sep=0.5pt]{$y_1\!\!-\!\!a_{2}$};
\draw(3,2)node[fill=white,inner sep=0.5pt]{$y_1\!\!-\!\!a_{4}$};
\draw(3,3)node[fill=white,inner sep=0.5pt]{$y_1\!\!-\!\!a_{5}$};
%\draw(4,-5)node[fill=white,inner sep=0.5pt]{$x_3\!\!+\!\!a_{\ov3}$};
\draw(4,-2)node[fill=white,inner sep=0.5pt]{$x_3\!\!+\!\!a_{0}$};
\draw(4,-1)node[fill=white,inner sep=0.5pt]{$x_3\!\!+\!\!a_{1}$};
\draw(4,1)node[fill=white,inner sep=0.5pt]{$x_3\!\!+\!\!a_{3}$};
\draw(5,3)node[fill=white,inner sep=0.5pt]{$y_2\!\!-\!\!a_{6}$};
\draw(5,4)node[fill=white,inner sep=0.5pt]{$y_2\!\!-\!\!a_{6}$};
\draw(6,-5)node[fill=white,inner sep=0.5pt]{$y_3\!\!-\!\!a_{\ov4}$};
\draw(6,0)node[fill=white,inner sep=0.5pt]{$y_3\!\!-\!\!a_{1}$};
\draw(7,-3)node[fill=white,inner sep=0.5pt]{$y_4\!\!-\!\!a_{\ov3}$};
\draw(8,-4)node[fill=white,inner sep=0.5pt]{$x_4\!\!+\!\!a_{\ov4}$};
\draw(8,-2)node[fill=white,inner sep=0.5pt]{$x_4\!\!+\!\!a_{\ov2}$};
\draw(8,-1)node[fill=white,inner sep=0.5pt]{$x_4\!\!+\!\!a_{\ov1}$};
\draw(8,5)node[fill=white,inner sep=0.5pt]{$x_4\!\!+\!\!a_{5}$};
\draw(9,-3)node[fill=white,inner sep=0.5pt]{$x_5\!\!+\!\!a_{\ov2}$};
\draw(9,2)node[fill=white,inner sep=0.5pt]{$x_5\!\!+\!\!a_{3}$};
%node{$k=$};
\draw(0,-7.5)node{$r$};\draw(0,-7.0)node{$\sigma(r)$};%\draw(0,-6.9)node{$x_k,y_\ell$};
\foreach \i in {1,2,4,8,9} \draw(\i,-7.5)node{$\blue{\i}$};
\foreach \i in {3,5,6,7} \draw(\i,-7.5)node{$\red{\i}$};
\foreach \i in {7,8} \draw(\i,-7.0)node{$0$};\foreach \i in {1,6,9} \draw(\i,-7.0)node{$1$};\foreach \i in {2,3,4,5} \draw(\i,-7.0)node{$2$};
\draw(1,-6.5)node{$\blue{x_1}$};\draw(2,-6.5)node{$\blue{x_2}$};\draw(4,-6.5)node{$\blue{x_3}$};\draw(8,-6.5)node{$\blue{x_4}$};\draw(9,-6.5)node{$\blue{x_5}$};
\draw(3,-6.5)node{$\red{y_1}$};\draw(5,-6.5)node{$\red{y_2}$};\draw(6,-6.5)node{$\red{y_3}$};\draw(7,-6.5)node{$\red{y_4}$};
%node{$c=$}
\draw(-0.9,-7)node{$c=$};
\foreach \j in {2,3,5,6} \draw(-0.9,-\j)node{$\green{\ov\j}$};
\foreach \j in {1,4} \draw(-0.9,-\j)node{$\brown{\ov\j}$};
\foreach \j in {0,4} \draw(-0.9,\j)node{$\green{\j}$};
\foreach \j in {1,2,3,5} \draw(-0.9,\j)node{$\brown{\j}$};
\draw(10.9,-7)node{$c=$};
\foreach \j in {2,3,5,6} \draw(10.9,-\j)node{$\green{\ov\j}$};
\foreach \j in {1,4} \draw(10.9,-\j)node{$\brown{\ov\j}$};
\foreach \j in {0,4} \draw(10.9,\j)node{$\green{\j}$};
\foreach \j in {1,2,3,5} \draw(10.9,\j)node{$\brown{\j}$};
%node{$a_m$}
%\draw(-0.9,-7)node{$c=$};
\foreach \j in {1,...,5} \draw(-0.4,-\j-1)node{$a_{\ov\j}$};
\foreach \j in {0,...,6} \draw(-0.4,\j-1)node{$a_{\j}$};
%\draw(10.3,-7)node{$c=$};
\foreach \j in {1,...,5} \draw(10.4,-\j-1)node{$a_{\ov\j}$};
\foreach \j in {0,...,6} \draw(10.4,\j-1)node{$a_{\j}$};
\end{tikzpicture}
}}
\end{equation}

With this Definition~\ref{Def-fact-sfnxy} we have
\begin{Proposition}\label{Prop-fact-ssfn}
The factorial skew Schur function $s^{R'}_{\lambda/\mu}(\x,\y|\a)$ is 
(i) independent of $R'$, ie. independent of the choice of $M$ and $N$; 
(ii) invariant under permutations of $(x_1,x_2,\ldots,x_m)$;
(iii) invariant under permutations of $(y_1,y_2,\ldots,y_n)$;
(iv) independent of $t$ if $x_k=t$ and $y_\ell=-t$ for any $k\in\{1,2,\ldots,m\}$ and any $\ell\in\{1,2,\ldots,n\}$.
\end{Proposition}

\noindent{\bf Proof}:
We shall consider consecutive letters $r$ and $r+1$ in the alphabet $R$
corresponding to a pair of consecutive horizontal lines in the lattice.
For typographical convenience we shall set $s=r+1$. 

There are four cases to consider, depending on whether $r$ and $s$ are either unprimed or primed in $R'$. 
In each case the boxes of any $T\in {\cal T}_{R'}^{\lambda/\mu}$ containing entries of the appropriate
types form a skew shape, a typical connected segment of which takes one of the following four forms.
\begin{equation}\label{eqn-tab-4}
%horizontal strip r<s
\begin{array}{cccc}
\vcenter{\hbox{
\begin{tikzpicture}[x={(0in,-0.2in)},y={(0.18in,0in)}] % matrix coordinate
%shape 
\foreach \j in {8,...,8} \draw[thick] (1,\j) rectangle +(-1,-1); 
\foreach \j in {1,...,8} \draw[thick] (2,\j) rectangle +(-1,-1);
\foreach \j in {1,...,2} \draw[thick] (3,\j) rectangle +(-1,-1);
%tableau entries
\foreach \j in {8,...,8} \draw(1-0.5,\j-0.5) node{$\blue{r}$};
\foreach \j in {1,...,5} \draw(2-0.5,\j-0.5) node{$\blue{r}$};
\foreach \j in {6,...,8} \draw(2-0.5,\j-0.5) node{$\red{s}$};
\foreach \j in {1,...,2}\draw(3-0.5,\j-0.5) node{$\red{s}$};
\end{tikzpicture}
}}
&
%vertical strip r'<s'
\vcenter{\hbox{
\begin{tikzpicture}[x={(0.18in,0in)},y={(0in,-0.2in)}] % matrix coordinate
%shape 
\foreach \j in {8,...,8} \draw[thick] (1,\j) rectangle +(-1,-1); 
\foreach \j in {1,...,8} \draw[thick] (2,\j) rectangle +(-1,-1);
\foreach \j in {1,...,2} \draw[thick] (3,\j) rectangle +(-1,-1);
%tableau entries
\foreach \j in {8,...,8} \draw(1-0.5,\j-0.5) node{$\blue{r'}$};
\foreach \j in {1,...,5} \draw(2-0.5,\j-0.5) node{$\blue{r'}$};
\foreach \j in {6,...,8} \draw(2-0.5,\j-0.5) node{$\red{s'}$};
\foreach \j in {1,...,2}\draw(3-0.5,\j-0.5) node{$\red{s'}$};
\end{tikzpicture}
}}
&
%strip r<s'
\vcenter{\hbox{
\begin{tikzpicture}[x={(0in,-0.2in)},y={(0.18in,0in)}] % matrix coordinate
%shape lambda
\foreach \j in {5,...,6} \draw[thick] (1,\j) rectangle +(-1,-1); 
\foreach \j in {5,...,5} \draw[thick] (2,\j) rectangle +(-1,-1);
\foreach \j in {5,...,5} \draw[thick] (3,\j) rectangle +(-1,-1);
\foreach \j in {5,...,5} \draw[thick] (4,\j) rectangle +(-1,-1);
\foreach \j in {2,...,5} \draw[thick] (5,\j) rectangle +(-1,-1);
\foreach \j in {2,...,2} \draw[thick] (6,\j) rectangle +(-1,-1);
%\foreach \j in {1,...,2} \draw[thick] (7,\j) rectangle +(-1,-1);
%tableau entries
\draw(1-0.5,5-0.5) node{$\blue{r}$};\draw(1-0.5,6-0.5) node{$\red{s'}$};
\draw(2-0.5,5-0.5) node{$\red{s'}$};
\draw(3-0.5,5-0.5) node{$\red{s'}$};
\draw(4-0.5,5-0.5) node{$\red{s'}$};
\draw(5-0.5,2-0.5) node{$\blue{r}$};\draw(5-0.5,3-0.5) node{$\blue{r}$};\draw(5-0.5,4-0.5) node{$\blue{r}$};\draw(5-0.5,5-0.5) node{$\red{s'}$};
\draw(6-0.5,2-0.5) node{$\red{s'}$};
%\draw(7-0.5,1-0.5) node{$\blue{r}$};\draw(7-0.5,2-0.5) node{$\red{s'}$};
\end{tikzpicture}
}}
&
%strip r'<s
\vcenter{\hbox{
\begin{tikzpicture}[x={(0in,-0.2in)},y={(0.18in,0in)}] % matrix coordinate
%shape lambda
\foreach \j in {5,...,6} \draw[thick] (1,\j) rectangle +(-1,-1); 
\foreach \j in {5,...,5} \draw[thick] (2,\j) rectangle +(-1,-1);
\foreach \j in {5,...,5} \draw[thick] (3,\j) rectangle +(-1,-1);
\foreach \j in {5,...,5} \draw[thick] (4,\j) rectangle +(-1,-1);
\foreach \j in {2,...,5} \draw[thick] (5,\j) rectangle +(-1,-1);
\foreach \j in {2,...,2} \draw[thick] (6,\j) rectangle +(-1,-1);
%\foreach \j in {1,...,2} \draw[thick] (7,\j) rectangle +(-1,-1);
%tableau entries
\draw(1-0.5,5-0.5) node{$\blue{r'}$};\draw(1-0.5,6-0.5) node{$\red{s}$};
\draw(2-0.5,5-0.5) node{$\blue{r'}$};
\draw(3-0.5,5-0.5) node{$\blue{r'}$};
\draw(4-0.5,5-0.5) node{$\blue{r'}$};
\draw(5-0.5,2-0.5) node{$\blue{r'}$};\draw(5-0.5,3-0.5) node{$\red{s}$};\draw(5-0.5,4-0.5) node{$\red{s}$};\draw(5-0.5,5-0.5) node{$\red{s}$};
\draw(6-0.5,2-0.5) node{$\blue{r'}$};
%\draw(7-0.5,1-0.5) node{$\blue{r}$};\draw(7-0.5,2-0.5) node{$\red{s'}$};
\end{tikzpicture}
}}
\end{array}
\end{equation}
Each of these supertableaux can be mapped bijectively to non-intersecting paths in a lattice path model based on an outside decomposition. 
It is convenient to choose outside decompositions based on the 
following cutting strips:
\begin{equation}\label{eqn-cuttingstrips-4}
\begin{array}{cccc}
\vcenter{\hbox{
\begin{tikzpicture}[x={(0in,-0.2in)},y={(0.18in,0in)}] % matrix coordinate
%shape 
\foreach \j in {-2,...,7} \draw[thick] (1,\j) rectangle +(-1,-1); 
%content
\foreach \j in {0,...,7} \draw(1-0.5,\j-0.5) node{$\brown{\j}$};
\foreach \j in {1,...,2} \draw(1-0.5,-\j-0.5) node{$\brown{\ov{\j}}$};
\end{tikzpicture}
}}
&
%vertical strip r'<s'
\vcenter{\hbox{
\begin{tikzpicture}[x={(0.18in,0in)},y={(0in,-0.2in)}] % matrix coordinate
%shape 
\foreach \j in {-2,...,7} \draw[thick] (1,\j) rectangle +(-1,-1); 
%content
\foreach \j in {0,...,7} \draw(1-0.5,\j-0.5) node{$\green{\ov{\j}}$};
\foreach \j in {1,...,2}\draw(1-0.5,-\j-0.5) node{$\green{\j}$};
\end{tikzpicture}
}}
&
%strip r<s'
\vcenter{\hbox{
\begin{tikzpicture}[x={(0in,-0.2in)},y={(0.18in,0in)}] % matrix coordinate
%shape lambda
\foreach \j in {5,...,6} \draw[thick] (1,\j) rectangle +(-1,-1); 
\foreach \j in {5,...,5} \draw[thick] (2,\j) rectangle +(-1,-1);
\foreach \j in {5,...,5} \draw[thick] (3,\j) rectangle +(-1,-1);
\foreach \j in {5,...,5} \draw[thick] (4,\j) rectangle +(-1,-1);
\foreach \j in {2,...,5} \draw[thick] (5,\j) rectangle +(-1,-1);
\foreach \j in {2,...,2} \draw[thick] (6,\j) rectangle +(-1,-1);
%\foreach \j in {1,...,2} \draw[thick] (7,\j) rectangle +(-1,-1);
%tableau entries
\draw(1-0.5,5-0.5) node{$\green{3}$};\draw(1-0.5,6-0.5) node{$\brown{4}$};
\draw(2-0.5,5-0.5) node{$\green{2}$};
\draw(3-0.5,5-0.5) node{$\green{1}$};
\draw(4-0.5,5-0.5) node{$\green{0}$};
\draw(5-0.5,2-0.5) node{$\green{\ov4}$};\draw(5-0.5,3-0.5) node{$\brown{\ov3}$};\draw(5-0.5,4-0.5) node{$\brown{\ov2}$};\draw(5-0.5,5-0.5) node{$\brown{\ov1}$};
\draw(6-0.5,2-0.5) node{$\green{\ov5}$};
%\draw(7-0.5,1-0.5) node{$\blue{r}$};\draw(7-0.5,2-0.5) node{$\red{s'}$};
\end{tikzpicture}
}}
&
%strip r'<s
\vcenter{\hbox{
\begin{tikzpicture}[x={(0in,-0.2in)},y={(0.18in,0in)}] % matrix coordinate
%shape lambda
\foreach \j in {5,...,6} \draw[thick] (1,\j) rectangle +(-1,-1); 
\foreach \j in {5,...,5} \draw[thick] (2,\j) rectangle +(-1,-1);
\foreach \j in {5,...,5} \draw[thick] (3,\j) rectangle +(-1,-1);
\foreach \j in {5,...,5} \draw[thick] (4,\j) rectangle +(-1,-1);
\foreach \j in {2,...,5} \draw[thick] (5,\j) rectangle +(-1,-1);
\foreach \j in {2,...,2} \draw[thick] (6,\j) rectangle +(-1,-1);
%\foreach \j in {1,...,2} \draw[thick] (7,\j) rectangle +(-1,-1);
%tableau entries
\draw(1-0.5,5-0.5) node{$\green{3}$};\draw(1-0.5,6-0.5) node{$\brown{4}$};
\draw(2-0.5,5-0.5) node{$\green{2}$};
\draw(3-0.5,5-0.5) node{$\green{1}$};
\draw(4-0.5,5-0.5) node{$\green{0}$};
\draw(5-0.5,2-0.5) node{$\green{\ov4}$};\draw(5-0.5,3-0.5) node{$\brown{\ov3}$};\draw(5-0.5,4-0.5) node{$\brown{\ov2}$};\draw(5-0.5,5-0.5) node{$\brown{\ov1}$};
\draw(6-0.5,2-0.5) node{$\green{\ov5}$};
%\draw(7-0.5,1-0.5) node{$\blue{r}$};\draw(7-0.5,2-0.5) node{$\red{s'}$};
\end{tikzpicture}
}}
\end{array}
\end{equation}
The precise content of each box in the skew shapes of (\ref{eqn-tab-4}) is %of course 
inherited from the supertableaux $T$ 
of which they form a part. This varies from case to case, with each variation obtainable from another through
the application of a suitable power of the shift operator $\tau$. For the purposes of illustration, the content
ascribed to the boxes of the cutting strips of (\ref{eqn-cuttingstrips-4}) is that which would follow from treating
each skew shaped tableau as a stand-alone skew supertableau in its own right. Then the four tableaux of (\ref{eqn-tab-4}) 
can be realised in the lattice framework as follows:

First the case $r$ and $s$ both unprimed in $R'$, for which Definition~\ref{Def-sigma} implies 
$\sigma(s)=\sigma(r)+1$.
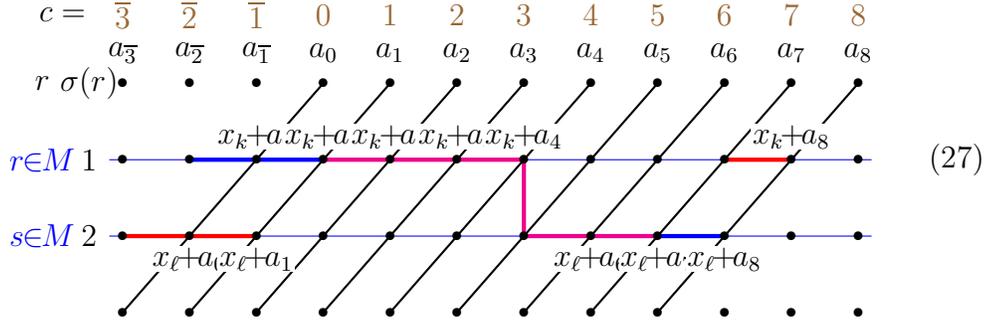
\begin{equation}\label{eqn-H}
\vcenter{\hbox{
\begin{tikzpicture}[x={(0in,-0.4in)},y={(0.35in,0in)}]%[x={(0in,-0.42in)},y={(0.42in,0in)}] % matrix coordinates
%horizontal x_k and y_l 
\foreach \i in {1,2} \draw[draw=blue](\i,-3.2)--(\i,8.2);
%paths of constant a_m
\foreach \j in {0,...,8} \draw[thick](0,\j)--(3,\j-3);
%path
\draw[draw=red,ultra thick] (1,6)--(1,7);
\draw[draw=blue,ultra thick] (1,-2)--(1,0);
\draw[draw=magenta,ultra thick] (1,0)--(1,3)--(2,3)--(2,5);
\draw[draw=blue,ultra thick] (2,5)--(2,6);
\draw[draw=red,ultra thick] (2,-3)--(2,-1);
%lattice grid
\foreach \i in {0,...,3} \foreach \j in {-3,...,8} \draw(\i,\j)node{$\sc\bullet$};
%\foreach \i in {1,2} \draw(\i,-7)node{$\sc\bullet$};
%weights
\draw(1-0.3,7)node[fill=white,inner sep=0.5pt]{$x_k\!\!+\!\!a_{8}$};
\draw(1-0.3,-1)node[fill=white,inner sep=0.5pt]{$x_k\!\!+\!\!a_{0}$};
\draw(1-0.3,0)node[fill=white,inner sep=0.5pt]{$x_k\!\!+\!\!a_{1}$};
\draw(1-0.3,1)node[fill=white,inner sep=0.5pt]{$x_k\!\!+\!\!a_{2}$};
\draw(1-0.3,2)node[fill=white,inner sep=0.5pt]{$x_k\!\!+\!\!a_{3}$};
\draw(1-0.3,3)node[fill=white,inner sep=0.5pt]{$x_k\!\!+\!\!a_{4}$};
\draw(2+0.3,4)node[fill=white,inner sep=0.5pt]{$x_\ell\!\!+\!\!a_{6}$};
\draw(2+0.3,5)node[fill=white,inner sep=0.5pt]{$x_\ell\!\!+\!\!a_{7}$};
\draw(2+0.3,6)node[fill=white,inner sep=0.5pt]{$x_\ell\!\!+\!\!a_{8}$};
\draw(2+0.3,-2)node[fill=white,inner sep=0.5pt]{$x_\ell\!\!+\!\!a_{0}$};
\draw(2+0.3,-1)node[fill=white,inner sep=0.5pt]{$x_\ell\!\!+\!\!a_{1}$};
%%
%node{$k=$};
\draw(0,-4.2)node{$r$};\draw(0,-3.5)node{$\sigma(r)$};%\draw(0,-6.9)node{$x_k,y_\ell$};
\draw(1,-4.2)node{$\blue{r\!\!\in\!\!M}$};\draw(1,-3.5)node{$1$};
\draw(2,-4.2)node{$\blue{s\!\!\in\!\!M}$};\draw(2,-3.5)node{$2$};
%node{$c=$}
\draw(-0.9,-3.9)node{$c=$};
\foreach \j in {1,...,3} \draw(-0.9,-\j)node{$\brown{\ov\j}$};
\foreach \j in {0,...,8} \draw(-0.9,\j)node{$\brown{\j}$};
%node{$a_m$}
%\draw(-0.9,-7)node{$c=$};
\foreach \j in {1,2,3} \draw(-0.4,-\j)node{$a_{\ov\j}$};
\foreach \j in {0,...,8} \draw(-0.4,\j)node{$a_{\j}$};
\end{tikzpicture}
}}
\end{equation}
Here it is important to note that with $r<s=r+1$ it is the condition $\sigma(s)=\sigma(r)+1$ that leads
to the lines of constant $a_t$ sloping diagonally down to the left. Furthermore it is this condition 
that guarantees that the sum of the weights of the boxes of each vertical pair with entries $r$ above $s$,
represented here by blue and red edges on the same diagonal, here labelled $a_t$ for $t=0$, $1$ and $8$, takes the form
$(x_k+a_t)(x_\ell+a_t)$ that is manifestly symmetric in $x_k$ and $x_\ell$. 

To be sure that $s_{\lambda/\mu}(\x,\y|\a)$ shares this symmetry it is necessary to go further. 
We refer to those entries that are not in a vertical pair as being unpaired. 
They constitute in general a number of horizontal strips.
In the non-factorial case symmetry is ensured by the observation~\cite{BK,KP}, that interchanging the 
number $\#(r)$ of unpaired entries $r$ of weight $x_k$ and the number $\#(s)$ of unpaired entries
of weight $x_\ell$ in every row of a supertableau $T$ yields, if these numbers differ, a supertableau 
$\tilde{T}$ such that the sum of the weights of $T$ and $\tilde{T}$ is symmetric in $x_k$ and $x_\ell$. 
If the numbers do not differ the weight of $T$ is of course symmetric.

The factorial case involves summing over all supertableaux $T$ for which $\#(r)+\#(s)$, the 
total number of entries $r$ and $s$ in any row, is fixed. In the example of (\ref{eqn-H}) it 
involves summing the contributions of a sequence of subtableaux or corresponding paths of the type:
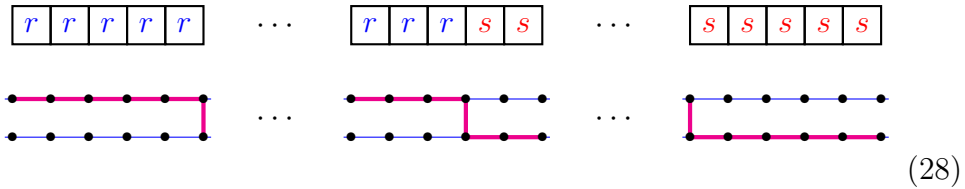
\begin{equation}
\begin{array}{ccccc}
\vcenter{\hbox{
\begin{tikzpicture}[x={(0in,-0.2in)},y={(0.2in,0in)}] % matrix coordinate
%shape 
\foreach \j in {1,...,5} \draw[thick] (2,\j) rectangle +(-1,-1);
%tableau entries
\foreach \j in {1,...,5} \draw(2-0.5,\j-0.5) node{$\blue{r}$};
%\foreach \j in {6,...,8} \draw(2-0.5,\j-0.5) node{$\red{s}$};
\end{tikzpicture}
}}
&\cdots
&
\vcenter{\hbox{
\begin{tikzpicture}[x={(0in,-0.2in)},y={(0.2in,0in)}] % matrix coordinate
%shape 
\foreach \j in {1,...,5} \draw[thick] (2,\j) rectangle +(-1,-1);
%tableau entries
\foreach \j in {1,...,3} \draw(2-0.5,\j-0.5) node{$\blue{r}$};
\foreach \j in {4,...,5} \draw(2-0.5,\j-0.5) node{$\red{s}$};
\end{tikzpicture}
}}
&\cdots
&
\vcenter{\hbox{
\begin{tikzpicture}[x={(0in,-0.2in)},y={(0.2in,0in)}] % matrix coordinate
%shape 
\foreach \j in {1,...,5} \draw[thick] (2,\j) rectangle +(-1,-1);
%tableau entries
%\foreach \j in {1,...,5} \draw(2-0.5,\j-0.5) node{$\blue{r}$};
\foreach \j in {1,...,5} \draw(2-0.5,\j-0.5) node{$\red{s}$};
\end{tikzpicture}
}}
\cr\cr
\vcenter{\hbox{
\begin{tikzpicture}[x={(0in,-0.2in)},y={(0.2in,0in)}] % matrix coordinates
%horizontal x_k and y_l 
\foreach \i in {1,2} \draw[draw=blue](\i,-0.2)--(\i,5.2);
%path
\draw[draw=magenta,ultra thick] (1,0)--(1,5)--(2,5);
%lattice grid
\foreach \i in {1,2} \foreach \j in {0,...,5} \draw(\i,\j)node{$\sc\bullet$};
\end{tikzpicture}
}}
&\cdots
&
\vcenter{\hbox{
\begin{tikzpicture}[x={(0in,-0.2in)},y={(0.2in,0in)}] % matrix coordinates
%horizontal x_k and y_l 
\foreach \i in {1,2} \draw[draw=blue](\i,-0.2)--(\i,5.2);
%path
\draw[draw=magenta,ultra thick] (1,0)--(1,3)--(2,3)--(2,5);
%lattice grid
\foreach \i in {1,2} \foreach \j in {0,...,5} \draw(\i,\j)node{$\sc\bullet$};
\end{tikzpicture}
}}
&\cdots
&
\vcenter{\hbox{
\begin{tikzpicture}[x={(0in,-0.2in)},y={(0.2in,0in)}] % matrix coordinates
%horizontal x_k and y_l 
\foreach \i in {1,2} \draw[draw=blue](\i,-0.2)--(\i,5.2);
%path
\draw[draw=magenta,ultra thick] (1,0)--(2,0)--(2,5);
%lattice grid
\foreach \i in {1,2} \foreach \j in {0,...,5} \draw(\i,\j)node{$\sc\bullet$};
\end{tikzpicture}
}}
\end{array}
\end{equation}
Each horizontal strip of length $m$ then contributes a factor 
$s_{(m)}(x_k,x_\ell,\0|\tau^p\a)$ for some shift parameter $p$. But this is just $s_{(m)}(x_k,x_\ell|\tau^p\a)$, which is 
an ordinary (not-supersymmetric) factorial Schur function that is well-known to be symmetric 
in $x_k$ and $x_\ell$~\cite{Mac92,GG}.

Similarly, if $r$ and $s$ are both primed in $R'$ then Definition~\ref{Def-sigma} implies 
$\sigma(s)=\sigma(r)-1$, and the lattice path picture typically takes the form:
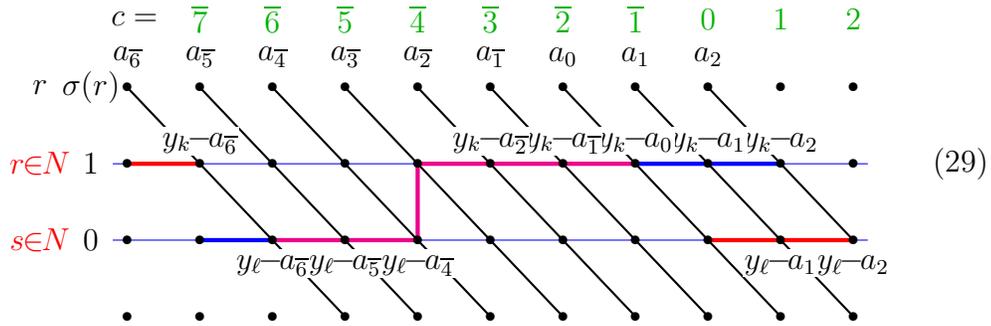
\begin{equation}\label{eqn-V}
\vcenter{\hbox{
\begin{tikzpicture}[x={(0in,-0.4in)},y={(0.38in,0in)}]%[x={(0in,-0.42in)},y={(0.42in,0in)}] % matrix coordinates
%horizontal x_k and y_l 
\foreach \i in {1,2} \draw[draw=blue](\i,-8.2)--(\i,2.2);
%paths of constant a_m
\foreach \j in {-8,...,-1} \draw[thick](0,\j)--(3,\j+3);
\draw[thick](0,0)--(2,2);
%path
\draw[draw=red,ultra thick] (2,0)--(2,2);
\draw[draw=blue,ultra thick] (1,-1)--(1,1);
\draw[draw=magenta,ultra thick] (2,-6)--(2,-4)--(1,-4)--(1,-1);
\draw[draw=red,ultra thick] (1,-8)--(1,-7);
\draw[draw=blue,ultra thick] (2,-7)--(2,-6);
%lattice grid
\foreach \i in {0,...,3} \foreach \j in {-8,...,2} \draw(\i,\j)node{$\sc\bullet$};
%\foreach \i in {1,2} \draw(\i,-7)node{$\sc\bullet$};
%weights
\draw(1-0.3,1)node[fill=white,inner sep=0.5pt]{$y_k\!\!-\!\!a_{2}$};
\draw(1-0.3,0)node[fill=white,inner sep=0.5pt]{$y_k\!\!-\!\!a_{1}$};
\draw(2+0.3,2)node[fill=white,inner sep=0.5pt]{$y_\ell\!\!-\!\!a_{2}$};
\draw(2+0.3,1)node[fill=white,inner sep=0.5pt]{$y_\ell\!\!-\!\!a_{1}$};
\draw(1-0.3,-1)node[fill=white,inner sep=0.5pt]{$y_k\!\!-\!\!a_{0}$};
\draw(1-0.3,-2)node[fill=white,inner sep=0.5pt]{$y_k\!\!-\!\!a_{\ov{1}}$};
\draw(1-0.3,-3)node[fill=white,inner sep=0.5pt]{$y_k\!\!-\!\!a_{\ov{2}}$};
\draw(2+0.3,-4)node[fill=white,inner sep=0.5pt]{$y_\ell\!\!-\!\!a_{\ov{4}}$};
\draw(2+0.3,-5)node[fill=white,inner sep=0.5pt]{$y_\ell\!\!-\!\!a_{\ov{5}}$};
\draw(2+0.3,-6)node[fill=white,inner sep=0.5pt]{$y_\ell\!\!-\!\!a_{\ov{6}}$};
\draw(1-0.3,-7)node[fill=white,inner sep=0.5pt]{$y_k\!\!-\!\!a_{\ov{6}}$};
%%
%node{$k=$};
\draw(0,-9.2)node{$r$};\draw(0,-8.5)node{$\sigma(r)$};%\draw(0,-6.9)node{$x_k,y_\ell$};
\draw(1,-9.2)node{$\red{r\!\!\in\!\!N}$};\draw(1,-8.5)node{$1$};
\draw(2,-9.2)node{$\red{s\!\!\in\!\!N}$};\draw(2,-8.5)node{$0$};
\draw(-0.9,-7.9)node{$c=$};
\foreach \j in {1,...,7} \draw(-0.9,-\j)node{$\green{\ov\j}$};
\foreach \j in {0,...,2} \draw(-0.9,\j)node{$\green{\j}$};
%node{$a_m$}
%\draw(-0.9,-7)node{$c=$};
\foreach \j in {1,...,6} \draw(-0.4,-\j-2)node{$a_{\ov\j}$};
\foreach \j in {0,...,2} \draw(-0.4,\j-2)node{$a_{\j}$};
\end{tikzpicture}
}}
\end{equation}
This time it is the condition $\sigma(s)=\sigma(r)-1$ that leads to the lines of constant $a_t$ sloping diagonally 
down to the right. This guarantees that the sum of the weights of the boxes of each horizontal pair with entry $s'$ 
to the right of an entry $r'$, represented here by blue and red edges on the same diagonal, labelled $a_t$ for $t=-6$, $1$ and $2$, takes the form
$(y_k-a_t)(y_\ell-a_t)$ that is manifestly symmetric in $y_k$ and $y_\ell$. 

To be sure that $s_{\lambda/\mu}(\x,\y|\a)$ shares this symmetry it is necessary to go further. 
As before we refer to those entries that are not in a horizontal pair as being unpaired. 
They constitute in general a number of vertical strips.
Once again in the non-factorial case symmetry is ensured by the same observation~\cite{BK,KP} as before
but now applied to the interchange of the numbers $\#(r')$ and $\#s')$ of unpaired primed entries $r'$ and $s'$,
respectively, in columns. Again as before, the factorial case involves summing over all 
supertableaux $T$ for which in any column $\#(r')+\#(s')$ is fixed. In the example of (\ref{eqn-V}) it 
involves summing the contributions of a sequence of subtableaux or corresponding paths of the type:
\begin{equation}
\begin{array}{ccccc}
\vcenter{\hbox{
\begin{tikzpicture}[y={(0in,-0.2in)},x={(+0.2in,0in)}] % matrix coordinate
%shape 
\foreach \j in {1,...,5} \draw[thick] (2,\j) rectangle +(-1,-1);
%tableau entries
\foreach \j in {1,...,5} \draw(2-0.5,\j-0.5) node{$\blue{r'}$};
%\foreach \j in {6,...,8} \draw(2-0.5,\j-0.5) node{$\red{s}$};
\end{tikzpicture}
}}
&\cdots
&
\vcenter{\hbox{
\begin{tikzpicture}[y={(0in,-0.2in)},x={(+0.2in,0in)}] % matrix coordinate
%shape 
\foreach \j in {1,...,5} \draw[thick] (2,\j) rectangle +(-1,-1);
%tableau entries
\foreach \j in {1,...,3} \draw(2-0.5,\j-0.5) node{$\blue{r'}$};
\foreach \j in {4,...,5} \draw(2-0.5,\j-0.5) node{$\red{s'}$};
\end{tikzpicture}
}}
&\cdots
&
\vcenter{\hbox{
\begin{tikzpicture}[y={(0in,-0.2in)},x={(0.2in,0in)}] % matrix coordinate
%shape 
\foreach \j in {1,...,5} \draw[thick] (2,\j) rectangle +(-1,-1);
%tableau entries
%\foreach \j in {1,...,5} \draw(2-0.5,\j-0.5) node{$\blue{r}$};
\foreach \j in {1,...,5} \draw(2-0.5,\j-0.5) node{$\red{s'}$};
\end{tikzpicture}
}}
\cr\cr
\vcenter{\hbox{
\begin{tikzpicture}[x={(0in,-0.2in)},y={(0.2in,0in)}] % matrix coordinates
%horizontal x_k and y_l 
\foreach \i in {1,2} \draw[draw=blue](\i,-0.2)--(\i,5.2);
%path
\draw[draw=magenta,ultra thick] (2,0)--(1,0)--(1,5);
%lattice grid
\foreach \i in {1,2} \foreach \j in {0,...,5} \draw(\i,\j)node{$\sc\bullet$};
\end{tikzpicture}
}}
&\cdots
&
\vcenter{\hbox{
\begin{tikzpicture}[x={(0in,-0.2in)},y={(0.2in,0in)}] % matrix coordinates
%horizontal x_k and y_l 
\foreach \i in {1,2} \draw[draw=blue](\i,-0.2)--(\i,5.2);
%path
\draw[draw=magenta,ultra thick] (2,0)--(2,2)--(1,2)--(1,5);
%lattice grid
\foreach \i in {1,2} \foreach \j in {0,...,5} \draw(\i,\j)node{$\sc\bullet$};
\end{tikzpicture}
}}
&\cdots
&
\vcenter{\hbox{
\begin{tikzpicture}[x={(0in,-0.2in)},y={(0.2in,0in)}] % matrix coordinates
%horizontal x_k and y_l 
\foreach \i in {1,2} \draw[draw=blue](\i,-0.2)--(\i,5.2);
%path
\draw[draw=magenta,ultra thick] (2,0)--(2,5)--(1,5);
%lattice grid
\foreach \i in {1,2} \foreach \j in {0,...,5} \draw(\i,\j)node{$\sc\bullet$};
\end{tikzpicture}
}}
\end{array}
\end{equation}
Each vertical strip of length $m$ contributes a factor $s_{(1^m)}(\0,y_k,y_\ell|\tau^p\a)$ for some shift parameter $p$. 
This can be seen to be identical to $s_{(m)}(y_\ell,y_k|-\tau^{p-m-1}\a)$,
but this is again an ordinary factorial Schur function necessarily symmetric in $y_k$ and $y_\ell$.

There remain the two cases with $r\in M,s\in N$ and $r\in N,s\in M$ as exemplified by the two
right-hand tableaux of (\ref{eqn-tab-4}), each in the form of a strip. For the left-hand strip
there exist just two tableaux of this shape with entries from the alphabet $\{r<s'\}$; the one illustrated
and one in which the entry $s'$ in the upper right-most box of maximum content is replaced by $r$.
The lattice path pictures are then as illustrated below where we have taken the liberty of superposing the pair 
of paths from the two allowed tableaux over the bulk of their length in magenta, while displaying the two 
different rightmost edges in blue and red.
\begin{equation}\label{eqn-MN}
\vcenter{\hbox{
\begin{tikzpicture}[x={(0in,-0.4in)},y={(0.38in,0in)}]%[x={(0in,-0.42in)},y={(0.42in,0in)}] % matrix coordinates
%horizontal x_k and y_l 
\draw[draw=blue](1,-6.2)--(1,4.2);
\draw[draw=red](2,-6.2)--(2,4.2);
%\foreach \i in {1} \foreach \j in {0.05} \fill[blue!40!white](\i-\j,-6.2)--(\i-\j,5.2)--(\i+\j,5.2)--(\i+\j,-6.2)--cycle;
%\foreach \i in {2} \foreach \j in {0.05} \fill[red!40!white](\i-\j,-6.2)--(\i-\j,5.2)--(\i+\j,5.2)--(\i+\j,-6.2)--cycle;
%paths of constant a_m
\foreach \j in {-5,...,4} \draw[thick](0,\j)--(3,\j);
%path
\draw[draw=magenta,ultra thick] (2,-6)--(2,-5)--(1,-4)--(1,-3)--(1,-2)--(2,-1)--(2,0)--(2,1)--(2,2)--(1,3);
\draw[draw=blue,ultra thick] (1,3)--(1,4);
\draw[draw=red,ultra thick] (1,3)--(2,4);
%\draw[draw=cyan,ultra thick] (0,0)--(2,0)--(2,1)--(3,2)--(4,2)--(5,3)--(5,4)--(8,4)--(8,5)--(10,5);
%\draw[draw=magenta,ultra thick] (0,2)--(2,2)--(3,3)--(0,3);
%lattice grid
\foreach \i in {0,...,3} \foreach \j in {-5,...,4} \draw(\i,\j)node{$\sc\bullet$};
\foreach \i in {1,2} \draw(\i,-6)node{$\sc\bullet$};
%weights
%\draw(1-0.3,-6)node[fill=white,inner sep=0.5pt]{$x_k\!\!+\!\!a_{\ov5}$};
%\draw(2+0.3,-5)node[fill=white,inner sep=0.5pt]{$y_\ell\!\!-\!\!a_{\ov4}$};
\draw(2+0.3,-5)node[fill=white,inner sep=0.5pt]{$y_\ell\!\!-\!\!a_{\ov4}$};
\draw(1-0.3,-4)node[fill=white,inner sep=0.5pt]{$x_k\!\!+\!\!a_{\ov3}$};
\draw(1-0.3,-3)node[fill=white,inner sep=0.5pt]{$x_k\!\!+\!\!a_{\ov2}$};
\draw(1-0.3,-2)node[fill=white,inner sep=0.5pt]{$x_k\!\!+\!\!a_{\ov1}$};
\draw(2+0.3,-1)node[fill=white,inner sep=0.5pt]{$y_\ell\!\!-\!\!a_{0}$};
\draw(2+0.3,0)node[fill=white,inner sep=0.5pt]{$y_\ell\!\!-\!\!a_{1}$};
\draw(2+0.3,1)node[fill=white,inner sep=0.5pt]{$y_\ell\!\!-\!\!a_{2}$};
\draw(2+0.3,2)node[fill=white,inner sep=0.5pt]{$y_\ell\!\!-\!\!a_{3}$};
\draw(1-0.3,3)node[fill=white,inner sep=0.5pt]{$x_k\!\!+\!\!a_{4}$};
\draw(1-0.3,4)node[fill=white,inner sep=0.5pt]{$x_k\!\!+\!\!a_{5}$};
\draw(2+0.3,4)node[fill=white,inner sep=0.5pt]{$y_k\!\!-\!\!a_{5}$};
%node{$k=$};
\draw(0,-7.2)node{$r$};\draw(0,-6.5)node{$\sigma(r)$};%\draw(0,-6.9)node{$x_k,y_\ell$};
\draw(1,-7.2)node{$\blue{r\!\!\in\!\!M}$};
\draw(2,-7.2)node{$\red{\s\!\!\in\!\!N}$};
\foreach \i in {1,2} \draw(\i,-6.5)node{$1$};
%\draw(1,-6.5)node{$\blue{x}$};
%\draw(2,-6.5)node{$\red{y}$};
%node{$c=$}
\draw(-0.9,-6)node{$c=$};
\foreach \j in {5,4} \draw(-0.9,-\j)node{$\green{\ov\j}$};
\foreach \j in {3,2,1} \draw(-0.9,-\j)node{$\brown{\ov\j}$};
\foreach \j in {0,1,2,3} \draw(-0.9,\j)node{$\green{\j}$};
\foreach \j in {4} \draw(-0.9,\j)node{$\brown{\j}$};
%node{$a_m$}
%\draw(-0.9,-7)node{$c=$};
\foreach \j in {1,...,4} \draw(-0.4,-\j-1)node{$a_{\ov\j}$};
\foreach \j in {0,...,5} \draw(-0.4,\j-1)node{$a_{\j}$};
\end{tikzpicture}
}}
\end{equation}
The sum of the weights of the blue and red edges is $(x_k+a_5)+(y_\ell-a_5)=x_k+y_\ell$,
so that the contribution to $s_{\lambda/\mu}(\x,\y|\a)$ of the only two tableaux of the
given strip shape with entries from the alphabet $\{r<s'\}$ contains a factor $(x_k+y_\ell)$ which
is clearly $0$ if $x_k=t=-\y_\ell$. Since these entries, if they appear at all, necessarily lie in
strips each sharing this property, it follows that $s_{\lambda/\mu}(\x,\y|\a)$ 
is independent of $t$ if $x_k=t=-y_\ell$, as required by the supersymmetry condition (iv).

A precisely analogous argument applies to the case $r\in N$ and $s\in M$.
The only significant difference is that this time the box of interest that may equally
well contain an entry $s'$ or $r$ lies at the lower lefthand end of the relevant strip
in the box of minimum content. The corresponding lattice picture is then as follows:
\begin{equation}\label{eqn-NM}
\vcenter{\hbox{
\begin{tikzpicture}[x={(0in,-0.4in)},y={(0.38in,0in)}] %[x={(0in,-0.42in)},y={(0.42in,0in)}] % matrix coordinates
%horizontal x_k and y_l 
\draw[draw=red](1,-6.2)--(1,4.2);
\draw[draw=blue](2,-6.2)--(2,4.2);
%\foreach \i in {1} \foreach \j in {0.05} \fill[blue!40!white](\i-\j,-6.2)--(\i-\j,5.2)--(\i+\j,5.2)--(\i+\j,-6.2)--cycle;
%\foreach \i in {2} \foreach \j in {0.05} \fill[red!40!white](\i-\j,-6.2)--(\i-\j,5.2)--(\i+\j,5.2)--(\i+\j,-6.2)--cycle;
%paths of constant a_m
\foreach \j in {-5,...,4} \draw[thick](0,\j)--(3,\j);
%path
\draw[draw=magenta,ultra thick] (1,-5)--(1,-4)--(2,-4)--(2,-3)--(2,-2)--(2,-1)--(1,-1)--(1,0)--(1,1)--(1,2)--(1,3)--(2,3)--(2,4);
\draw[draw=blue,ultra thick] (1,-6)--(1,-5);
\draw[draw=red,ultra thick] (2,-6)--(2,-5)--(1,-5);
%\draw[draw=cyan,ultra thick] (0,0)--(2,0)--(2,1)--(3,2)--(4,2)--(5,3)--(5,4)--(8,4)--(8,5)--(10,5);
%\draw[draw=magenta,ultra thick] (0,2)--(2,2)--(3,3)--(0,3);
%lattice grid
\foreach \i in {0,...,3} \foreach \j in {-5,...,4} \draw(\i,\j)node{$\sc\bullet$};
\foreach \i in {1,2} \draw(\i,-6)node{$\sc\bullet$};
%weights
%\draw(2+0.3,-6)node[fill=white,inner sep=0.5pt]{$x_k\!\!+\!\!a_{\ov6}$};
%\draw(1-0.3,-6)node[fill=white,inner sep=0.5pt]{$y_\ell\!\!-\!\!a_{\ov6}$};
\draw(2+0.3,-5)node[fill=white,inner sep=0.5pt]{$x_k\!\!+\!\!a_{\ov5}$};
\draw(1-0.3,-5)node[fill=white,inner sep=0.5pt]{$y_\ell\!\!-\!\!a_{\ov5}$};
\draw(1-0.3,-4)node[fill=white,inner sep=0.5pt]{$y_\ell\!\!-\!\!a_{\ov4}$};
\draw(2+0.3,-3)node[fill=white,inner sep=0.5pt]{$x_k\!\!+\!\!a_{\ov3}$};
\draw(2+0.3,-2)node[fill=white,inner sep=0.5pt]{$x_k\!\!+\!\!a_{\ov2}$};
\draw(2+0.3,-1)node[fill=white,inner sep=0.5pt]{$x_k\!\!+\!\!a_{\ov1}$};
\draw(1-0.3,0)node[fill=white,inner sep=0.5pt]{$y_\ell\!\!-\!\!a_{0}$};
\draw(1-0.3,1)node[fill=white,inner sep=0.5pt]{$y_\ell\!\!-\!\!a_{1}$};
\draw(1-0.3,2)node[fill=white,inner sep=0.5pt]{$y_\ell\!\!-\!\!a_{2}$};
\draw(1-0.3,3)node[fill=white,inner sep=0.5pt]{$y_\ell\!\!-\!\!a_{3}$};
\draw(2+0.3,4)node[fill=white,inner sep=0.5pt]{$x_k\!\!+\!\!a_{4}$};
%node{$k=$};
\draw(0,-7.2)node{$r$};\draw(0,-6.5)node{$\sigma(r)$};%\draw(0,-6.9)node{$x_k,y_\ell$};
\draw(1,-7.2)node{$\red{r\!\!\in\!\! N}$};
\draw(2,-7.2)node{$\blue{s\!\!\in\!\! M}$};
\foreach \i in {1,2} \draw(\i,-6.5)node{$0$};
%\draw(1,-6.5)node{$\blue{x}$};
%\draw(2,-6.5)node{$\red{y}$};
%node{$c=$}
\draw(-0.9,-6)node{$c=$};
\foreach \j in {5,4} \draw(-0.9,-\j)node{$\green{\ov\j}$};
\foreach \j in {3,2,1} \draw(-0.9,-\j)node{$\brown{\ov\j}$};
\foreach \j in {0,1,2,3} \draw(-0.9,\j)node{$\green{\j}$};
\foreach \j in {4} \draw(-0.9,\j)node{$\brown{\j}$};
%node{$a_m$}
%\draw(-0.9,-7)node{$c=$};
\foreach \j in {1,...,5} \draw(-0.4,-\j)node{$a_{\ov\j}$};
\foreach \j in {0,...,4} \draw(-0.4,\j)node{$a_{\j}$};
\end{tikzpicture}
}}
\end{equation}
The sum of the weights of the blue and red edges is $(x_k+a_{\ov5})+(y_\ell-a_{\ov5})=x_k+y_\ell$,
as required to show that the contribution to $s_{\lambda/\mu}(\x,\y|\a)$ of the only two tableaux of the
given strip shape with entries from the alphabet $\{r'<s\}$ is $0$ if $x_k=t=-y_\ell$. 
As before, since these entries, if they appear at all, necessarily lie in strips sharing this 
property, it follows once again that $s_{\lambda/\mu}(\x,\y|\a)$ 
is independent of $t$ if $x_k=t=-y_\ell$, as required by the supersymmetry condition (iv).

The map between these two pictures (\ref{eqn-MN}) and (\ref{eqn-NM}) is bijective and corresponds 
to the map from the pair of tableaux on the left below to the pair on the right:
\begin{equation}\label{eqn-cyclic}
\begin{array}{ccc}  %%ccccccc
\vcenter{\hbox{
\begin{tikzpicture}[x={(0in,-0.2in)},y={(0.15in,0in)}] %[x={(0in,-0.2in)},y={(0.2in,0in)}] % matrix coordinate
%shape lambda
\foreach \j in {5,...,6} \draw[thick] (1,\j) rectangle +(-1,-1); 
\foreach \j in {5,...,5} \draw[thick] (2,\j) rectangle +(-1,-1);
\foreach \j in {5,...,5} \draw[thick] (3,\j) rectangle +(-1,-1);
\foreach \j in {5,...,5} \draw[thick] (4,\j) rectangle +(-1,-1);
\foreach \j in {2,...,5} \draw[thick] (5,\j) rectangle +(-1,-1);
\foreach \j in {2,...,2} \draw[thick] (6,\j) rectangle +(-1,-1);
\foreach \j in {1,...,2} \draw[thick] (7,\j) rectangle +(-1,-1);
%tableau entries
\draw(1-0.5,5-0.5) node{$\blue{r}$};\draw(1-0.5,6-0.5) node{$\blue{r}$};
\draw(2-0.5,5-0.5) node{$\red{s'}$};
\draw(3-0.5,5-0.5) node{$\red{s'}$};
\draw(4-0.5,5-0.5) node{$\red{s'}$};
\draw(5-0.5,2-0.5) node{$\blue{r}$};\draw(5-0.5,3-0.5) node{$\blue{r}$};\draw(5-0.5,4-0.5) node{$\blue{r}$};\draw(5-0.5,5-0.5) node{$\red{s'}$};
\draw(6-0.5,2-0.5) node{$\red{s'}$};
\draw(7-0.5,1-0.5) node{$\blue{r}$};\draw(7-0.5,2-0.5) node{$\red{s'}$};
\end{tikzpicture}
}}
\!\!\!+\!\!\!  %&+&
\vcenter{\hbox{
\begin{tikzpicture}[x={(0in,-0.2in)},y={(0.15in,0in)}] %[x={(0in,-0.2in)},y={(0.2in,0in)}] % matrix coordinate
%shape lambda
\foreach \j in {5,...,6} \draw[thick] (1,\j) rectangle +(-1,-1); 
\foreach \j in {5,...,5} \draw[thick] (2,\j) rectangle +(-1,-1);
\foreach \j in {5,...,5} \draw[thick] (3,\j) rectangle +(-1,-1);
\foreach \j in {5,...,5} \draw[thick] (4,\j) rectangle +(-1,-1);
\foreach \j in {2,...,5} \draw[thick] (5,\j) rectangle +(-1,-1);
\foreach \j in {2,...,2} \draw[thick] (6,\j) rectangle +(-1,-1);
\foreach \j in {1,...,2} \draw[thick] (7,\j) rectangle +(-1,-1);
%tableau entries
\draw(1-0.5,5-0.5) node{$\blue{r}$};\draw(1-0.5,6-0.5) node{$\red{s'}$};
\draw(2-0.5,5-0.5) node{$\red{s'}$};
\draw(3-0.5,5-0.5) node{$\red{s'}$};
\draw(4-0.5,5-0.5) node{$\red{s'}$};
\draw(5-0.5,2-0.5) node{$\blue{r}$};\draw(5-0.5,3-0.5) node{$\blue{r}$};\draw(5-0.5,4-0.5) node{$\blue{r}$};\draw(5-0.5,5-0.5) node{$\red{s'}$};
\draw(6-0.5,2-0.5) node{$\red{s'}$};
\draw(7-0.5,1-0.5) node{$\blue{r}$};\draw(7-0.5,2-0.5) node{$\red{s'}$};
\end{tikzpicture}
}}
&\rightarrow&
\vcenter{\hbox{
\begin{tikzpicture}[x={(0in,-0.2in)},y={(0.15in,0in)}] %[x={(0in,-0.2in)},y={(0.2in,0in)}] % matrix coordinate
%shape lambda
\foreach \j in {5,...,6} \draw[thick] (1,\j) rectangle +(-1,-1); 
\foreach \j in {5,...,5} \draw[thick] (2,\j) rectangle +(-1,-1);
\foreach \j in {5,...,5} \draw[thick] (3,\j) rectangle +(-1,-1);
\foreach \j in {5,...,5} \draw[thick] (4,\j) rectangle +(-1,-1);
\foreach \j in {2,...,5} \draw[thick] (5,\j) rectangle +(-1,-1);
\foreach \j in {2,...,2} \draw[thick] (6,\j) rectangle +(-1,-1);
\foreach \j in {1,...,2} \draw[thick] (7,\j) rectangle +(-1,-1);
%tableau entries
\draw(1-0.5,5-0.5) node{$\red{s'}$};\draw(1-0.5,6-0.5) node{$\blue{r}$};
\draw(2-0.5,5-0.5) node{$\red{s'}$};
\draw(3-0.5,5-0.5) node{$\red{s'}$};
\draw(4-0.5,5-0.5) node{$\red{s'}$};
\draw(5-0.5,2-0.5) node{$\red{s'}$};\draw(5-0.5,3-0.5) node{$\blue{r}$};\draw(5-0.5,4-0.5) node{$\blue{r}$};\draw(5-0.5,5-0.5) node{$\blue{r}$};
\draw(6-0.5,2-0.5) node{$\red{s'}$};
\draw(7-0.5,1-0.5) node{$\blue{r}$};\draw(7-0.5,2-0.5) node{$\blue{r}$};
\end{tikzpicture}
}}
\!\!\!+\!\!\! %&+&
\vcenter{\hbox{
\begin{tikzpicture}[x={(0in,-0.2in)},y={(0.15in,0in)}]%[x={(0in,-0.2in)},y={(0.2in,0in)}] % matrix coordinate
%shape lambda
\foreach \j in {5,...,6} \draw[thick] (1,\j) rectangle +(-1,-1); 
\foreach \j in {5,...,5} \draw[thick] (2,\j) rectangle +(-1,-1);
\foreach \j in {5,...,5} \draw[thick] (3,\j) rectangle +(-1,-1);
\foreach \j in {5,...,5} \draw[thick] (4,\j) rectangle +(-1,-1);
\foreach \j in {2,...,5} \draw[thick] (5,\j) rectangle +(-1,-1);
\foreach \j in {2,...,2} \draw[thick] (6,\j) rectangle +(-1,-1);
\foreach \j in {1,...,2} \draw[thick] (7,\j) rectangle +(-1,-1);
%tableau entries
\draw(1-0.5,5-0.5) node{$\red{s'}$};\draw(1-0.5,6-0.5) node{$\blue{r}$};
\draw(2-0.5,5-0.5) node{$\red{s'}$};
\draw(3-0.5,5-0.5) node{$\red{s'}$};
\draw(4-0.5,5-0.5) node{$\red{s'}$};
\draw(5-0.5,2-0.5) node{$\red{s'}$};\draw(5-0.5,3-0.5) node{$\blue{r}$};\draw(5-0.5,4-0.5) node{$\blue{r}$};\draw(5-0.5,5-0.5) node{$\blue{r}$};
\draw(6-0.5,2-0.5) node{$\red{s'}$};
\draw(7-0.5,1-0.5) node{$\red{s'}$};\draw(7-0.5,2-0.5) node{$\blue{r}$};
\end{tikzpicture}
}}\cr
\end{array}
\end{equation}
The transformation from a pair of tableaux on the left with entries $r<s'$ to a pair on the right with entries $s'<r$ 
can be seen to be a weight preserving map that cyclically permutes the entries 
by moving each entry one step either up or to the right, with the top rightmost entry moved to the bottom leftmost position. 
The resulting pair of tableaux are standard with respect to the order $s'<r$. 
Since all other entries remain fixed it can be seen that this map applied to all the strips containing entries 
$r$ and $s'$ maps bijectively all $T'\in{\cal T}^{\lambda/\mu}_{R'}$ based on an alphabet $R'=(\cdots<r<s'<\cdots)$
to corresponding $T''\in{\cal T}^{\lambda/\mu}_{R''}$ based on the alphabet $R''=(\cdot<s'<r<\cdots)$ obtained from $R'$ by
interchanging $r$ and $s'$. It follows that $s_{\lambda/\mu}^{R'}(\x,\y|\a)=s_{\lambda/\mu}^{R''}(\x,\y|\a)$.
Similarly the inverse map can be used to bijectively map all $T''\in{\cal T}^{\lambda/\mu}_{R''}$ to
all $T'\in{\cal T}^{\lambda/\mu}_{R'}$. 

For any choice of alphabet $R'$ involving both unprimed and primed there will necessarily exist at least one pair of neighbouring 
elements $r$ and $s'$ with either $r<s'$ or $s'<r$, and correspondingly there exists $k$ and $\ell$ such that the
supersymmetry condition (iv) is satisfied. Moreover, iterating the cyclic transformation exemplified by (\ref{eqn-cyclic}),
or its inverse, it is easy to see that the skew Schur function $s_{\lambda/\mu}^{R'}(\x,\y|\a)$ defined in terms of 
supertableaux based on an alphabet $R'=M\dot\cup N'$ for some fixed but arbitrary choice of $M$ and $N$ is identical to one 
based on the alphabet $\{1'<2'<\cdots<n'<n+1<n+2<\cdots<n+m\}$ or $\{1<2<\cdots<m<(m+1)'<(m+2)'<\cdots<(m+n)'\}$. 
It then follows from our symmetry properties in the $r<s$ and $r'<s'$ cases that both of the required supersymmetry conditions 
(ii) and (iii) are satisfied, while at the same time we can conclude that $s_{\lambda/\mu}^{R'}(\x,\y|\a)$ 
is independent of the choice of $M$ and $N$, so that condition (i) is also satisfied, and the 
superscript $R'$ is redundant since the the identification of which $N$ elements are primed is immaterial.
\qed

\section{Conclusion}\label{Sec-conclusion}

In the discussion of the factorial supersymmetric skew Schur functions 
of the previous section, it should be noted that it is the simple requirement of symmetry of the weights of
vertical pairs of entries $r$ and $s$, and of horizontal pairs of entries $r'$ and $s'$ that are
responsible for the key feature of the Definition~\ref{Def-sigma} whereby $\sigma(r)$ increases by $1$
and decreases by $1$ in the passage from $r$ to $s$, and from $r'$ to $s'$, respectively. The requirement 
of the supersymmetric condition (iv) on supertableaux containing strips of length $1$ having entries $r$
and $s'$, or $r'$ and $s$ is sufficient to show that $\sigma(r)$ must remain unchanged in the passage
from $r$ to $s'$ and from $r'$ to $s$. These observations are sufficient to fix Definition~\ref{Def-sigma}
completely.

What is perhaps remarkable is that having fixed this definition based on pairs or singletons of entries, it is
found that not only are the supersymmetry conditions (ii), (iii) and (iv) found to hold quite generally, but also
condition (i) regarding the independence of the resulting factorial supersymmetric skew Schur functions on the
choice of primed alphabet $R'$.

In particular this accommodates both the alphabet $\{1<2<\cdots<m<(m+1)'<(m+2)'<\cdots<(m+n)'\}$ used, 
with $(m+\ell)'$ replaced by $\ell'$, by Berele and Regev~\cite{BR87} in defining non-factorial 
supersymmetric Schur functions for the first time, and the alphabet $\{1'<2'<\cdots<n'<n+1<n+2<\cdots<n+m\}$,
established, with $n+k$ replaced by $k$, as an alternative by Remmel~\cite{Rem1,Rem2} in this same non-factorial case.
This was done by means of a variant of Sch\"utzenberger's {\em jeu de taquin}~\cite{Sch}
that moves unprimed entries individually step by step past primed entries towards the outer rim of 
the tableau. Such steps are weight preserving in the non-factorial case but not in the factorial case,
so cannot be used here. Instead, this has been accomplished here for factorial supersymmetric skew Schur functions 
through the iterative use of our cyclic map, which incidentally requires significantly fewer steps than
using {\em jeu de taquin}.

Even in the factorial case it should also be pointed out that the ordering of the elements in $M=\{i_1<i_2<\cdots<i_m\}$ and 
$N=\{j_i<j_2<\cdots<j_n\}$ does not necessarily require that $i_k$ and $j_\ell$ carry weights 
$(x_k+a_{\sigma(r)+c})$ and $(y_\ell-a_{\sigma(s)+c+1})$, respectively, as has been assumed here.
Indeed, it is clear from the supersymmetry conditions that the components of $\x$ and $\y$ may be independently
permuted in any way that one wishes. One particularly notable choice, already encountered in the Introduction, 
is that of Molev~\cite{Mol} which amounts to using the alphabet $\{n'<\cdots<2'<1'<1<2<\cdots<m\}$. In our formalism
this is equivalent to using the alphabet $R'=\{1'<2'<\cdots<n'<n+1<n+2<\cdots<n+m\}$ with entries 
carrying weights
\begin{equation}
\wgt(t_{ij})=\begin{cases}
                  x_k+a_{k+c-n}&\mbox{if $t_{ij}=n+k\in M$};\cr
									y_\ell-a_{\ell+c-n}&\mbox{if $t_{ij}=(n-\ell+1)'\in N'$}.\cr
             \end{cases}
\end{equation} 
Comparing this with the formula given by Molev~\cite{Mol} $\{n'<\cdots<2'<1'<1<2<\cdots<m\}$ for
which
\begin{equation}
\wgt(t_{ij})=\begin{cases}
                  x_k+a_{k+c}&\mbox{if $t_{ij}=k\in M$};\cr
									y_\ell-a_{\ell+c}&\mbox{if $t_{ij}=\ell'\in N'$},\cr
             \end{cases}
\end{equation} 
one can see that
\begin{equation}
  s_{\lambda/\mu}^{R'}(\x,\y|\a)=s_{\lambda/\mu}^{\Mol}(\x,\y|\tau^{-n}\a)\,,
\end{equation}
where the superscript $\Mol$ indicates the factorial supersymmetric skew Schur function of~\cite{Mol}
and the action of $\tau^{-n}$ on any $a_{s}$ in $\a$ gives $a_{s-n}$. Since our factorial parameters are 
arbitrary, such a fixed translation in the indices of $\a$ is immaterial provided one sticks to the same 
translation for all Schur functions in all identities such as those of our determinantal formulae 
based on outside decompositions.\

It follows by setting $\a=\0=(\ldots,0,0,\ldots)$ in Definition~\ref{Def-fact-sfnxy} that $s_{\lambda/\mu}^{R'}(\x,\y|\0)$
is independent of the choice of $M$ and $N$ thanks to Proposition~\ref{Prop-fact-ssfn}, and by comparing weights coincides with the classical
supersymmetric skew Schur function $s_{\lambda/\mu}(\x,\y)$ of (\ref{eqn-sfnxy}) regardless of this choice.

Returning to the tableau-based ninth variation generalised skew Schur functions, $s_{\lambda/\mu}^{R'}(\X,\Y)$, of Definition~\ref{Def-sfnXY-R'},
these can be used to recover as a special case a skew extension of the Schur version of objects introduced by Bachmann~\cite{Bac} that
are based on what he called ordered Young tableaux (OYTs). An ordered Young tableau is one in which entries from an alphabet $\{1<2<\cdots<n\}$
are weakly increasing across rows from left to right and down columns from top to bottom, and are strictly increasing down diagonals
from top left to bottom right. Let ${\cal OYT}^{\lambda/\mu}$ be the set of all such tableaux, $T$, of skew shape $F^{\lambda/\mu}$
and let the number of adjacent horizontal and vertical pairs of identical entries be $h(T)$ and $v(T)$, respectively. Then for
any parameters $F=(f_{kc})$ one can define the following skew version of Bachmann's Schur functions:
\begin{equation}\label{eqn-sfn-Ft}
   s_{\lambda/\mu}^{Bac}(F;t) = \sum_{T\in{\cal OTY}^{\lambda/\mu}} \ t^{v(T)}\,(1-t)^{h(T)} \prod_{(i,j)\in F^{\lambda/\mu}} \wgt(t_{ij})
\end{equation}
where if $t_{ij}=k$ then $\wgt(t_{ij})=f_{kc}$ with $c=j-i$. 

In our Definition~\ref{Def-sfnXY-R'} if we set $m=n$, $R'=\{1'<1<2<2'<\cdots<n'<n\}$, $x_{kc}=(1-t)f_{kc}$ and $y_{kc}=tf_{kc}$.
then we have
\begin{equation}\label{eqn-XY-Ft}
   s_{\lambda/\mu}^{R'}(\X,\Y)=s_{\lambda/\mu}^{\Bac}(F;t)\,. 
\end{equation}

To establish this one first notes that for any tableau $T\in{\cal OTY}^{\lambda/\mu}$ containing entries of fixed value, say $k$, they must
lie in a strip or a sequence of disconnected strips which we refer to as $k$-strips. For each constituent $k$-strip there is a $1-2$ map of entries $k$ to
corresponding strips of entries from $\{k'<k\}$ in a supertableau. This is exemplified by
\begin{equation}\label{eqn-OYTstrip}
\begin{array}{ccccccc} 
\vcenter{\hbox{
\begin{tikzpicture}[x={(0in,-0.2in)},y={(0.2in,0in)}] % matrix coordinate
%shape lambda
\foreach \j in {5,...,6} \draw[thick] (1,\j) rectangle +(-1,-1); 
\foreach \j in {5,...,5} \draw[thick] (2,\j) rectangle +(-1,-1);
\foreach \j in {5,...,5} \draw[thick] (3,\j) rectangle +(-1,-1);
\foreach \j in {5,...,5} \draw[thick] (4,\j) rectangle +(-1,-1);
\foreach \j in {2,...,5} \draw[thick] (5,\j) rectangle +(-1,-1);
\foreach \j in {2,...,2} \draw[thick] (6,\j) rectangle +(-1,-1);
\foreach \j in {1,...,2} \draw[thick] (7,\j) rectangle +(-1,-1);
%tableau entries
\draw(1-0.5,5-0.5) node{${k}$};\draw(1-0.5,6-0.5) node{${k}$};
\draw(2-0.5,5-0.5) node{${k}$};
\draw(3-0.5,5-0.5) node{${k}$};
\draw(4-0.5,5-0.5) node{${k}$};
\draw(5-0.5,2-0.5) node{${k}$};\draw(5-0.5,3-0.5) node{${k}$};\draw(5-0.5,4-0.5) node{${k}$};\draw(5-0.5,5-0.5) node{${k}$};
\draw(6-0.5,2-0.5) node{${k}$};
\draw(7-0.5,1-0.5) node{${k}$};\draw(7-0.5,2-0.5) node{${k}$};
\end{tikzpicture}
}}
&\mapsto&
\vcenter{\hbox{
\begin{tikzpicture}[x={(0in,-0.2in)},y={(0.2in,0in)}] % matrix coordinate
%shape lambda
\foreach \j in {5,...,6} \draw[thick] (1,\j) rectangle +(-1,-1); 
\foreach \j in {5,...,5} \draw[thick] (2,\j) rectangle +(-1,-1);
\foreach \j in {5,...,5} \draw[thick] (3,\j) rectangle +(-1,-1);
\foreach \j in {5,...,5} \draw[thick] (4,\j) rectangle +(-1,-1);
\foreach \j in {2,...,5} \draw[thick] (5,\j) rectangle +(-1,-1);
\foreach \j in {2,...,2} \draw[thick] (6,\j) rectangle +(-1,-1);
\foreach \j in {1,...,2} \draw[thick] (7,\j) rectangle +(-1,-1);
%tableau entries
\draw(1-0.5,5-0.5) node{$\red{k'}$};\draw(1-0.5,6-0.5) node{$\blue{k}$};
\draw(2-0.5,5-0.5) node{$\red{k'}$};
\draw(3-0.5,5-0.5) node{$\red{k'}$};
\draw(4-0.5,5-0.5) node{$\red{k'}$};
\draw(5-0.5,2-0.5) node{$\red{k'}$};\draw(5-0.5,3-0.5) node{$\blue{k}$};\draw(5-0.5,4-0.5) node{$\blue{k}$};\draw(5-0.5,5-0.5) node{$\blue{k}$};
\draw(6-0.5,2-0.5) node{$\red{k'}$};
\draw(7-0.5,1-0.5) node{$\blue{k}$};\draw(7-0.5,2-0.5) node{$\blue{k}$};
\end{tikzpicture}
}}
&+&
\vcenter{\hbox{
\begin{tikzpicture}[x={(0in,-0.2in)},y={(0.2in,0in)}] % matrix coordinate
%shape lambda
\foreach \j in {5,...,6} \draw[thick] (1,\j) rectangle +(-1,-1); 
\foreach \j in {5,...,5} \draw[thick] (2,\j) rectangle +(-1,-1);
\foreach \j in {5,...,5} \draw[thick] (3,\j) rectangle +(-1,-1);
\foreach \j in {5,...,5} \draw[thick] (4,\j) rectangle +(-1,-1);
\foreach \j in {2,...,5} \draw[thick] (5,\j) rectangle +(-1,-1);
\foreach \j in {2,...,2} \draw[thick] (6,\j) rectangle +(-1,-1);
\foreach \j in {1,...,2} \draw[thick] (7,\j) rectangle +(-1,-1);
%tableau entries
\draw(1-0.5,5-0.5) node{$\red{k'}$};\draw(1-0.5,6-0.5) node{$\blue{k}$};
\draw(2-0.5,5-0.5) node{$\red{k'}$};
\draw(3-0.5,5-0.5) node{$\red{k'}$};
\draw(4-0.5,5-0.5) node{$\red{k'}$};
\draw(5-0.5,2-0.5) node{$\red{k'}$};\draw(5-0.5,3-0.5) node{$\blue{k}$};\draw(5-0.5,4-0.5) node{$\blue{k}$};\draw(5-0.5,5-0.5) node{$\blue{k}$};
\draw(6-0.5,2-0.5) node{$\red{k'}$};
\draw(7-0.5,1-0.5) node{$\red{k'}$};\draw(7-0.5,2-0.5) node{$\blue{k}$};
\end{tikzpicture}
}}
\end{array}
\end{equation}
Apart from common factors $f_{kc}$ in the contribution to the weight, one for each box, the $t$-dependence on the left is $t^6(1-t)^5$ since 
for this subtableau $T$ we have $h(T)=5$ and $v(t)=6$, while the $t$-dependence on the right is $t^6(1-t)^6+t^7(1-t)^5=t^6(1-t)^5$ since the numbers of unprimed and primed entries on the right are both $6$ for the first subtableau and $5$ and $7$ for the second. This equality of weights comes about quite generally because for each such strip the lowest leftmost box of the strip of content $c$ may contain an entry $k$ which may be unprimed or uprimed on the right, leading to a combined
contribution to the weight of $(1-t)f_{kc}+tf_{kc}=f_{kc}$. At the same time, the rightmost and topmost entry $k$ of each horizontal and vertical pair of entries $k$ is necessarily mapped on the right to an unprimed and primed entry, $k$ and $k'$ respectively, of weights $(1-t)f_{kc}$ and $tf_{kc}$. 
Extracting all common factors of $t$ and $(1-t)$
for all possible $k$-strips for $k=1,2,\ldots,n$ then leads directly to the required weighting of $T\in{\cal OYT}^{\lambda/\mu}$, 
thereby establishing the identity (\ref{eqn-XY-Ft}). An immediate consequence of this identification of $s_{\lambda/\mu}^{\Bac}(F;t)$ 
with a special case of $s_{\lambda/\mu}^{R'}(\X,\Y)$ is that $s_{\lambda/\mu}^{\Bac}(F;t)$ must satisfy all the determinantal
identities of Theorem~\ref{The-HGXY} based on arbitrary outside decompositions, not just the Jacobi-Trudi and dual Jacobi Trudi identities 
of Theorem~5.1 in~\cite{Bac}.

In this same case $m=n$ with $R'=\{1'<1<2<2'<\cdots<n'<n\}$, if we revert to the sixth variation then we find 
from Definition~\ref{Def-sigma} that $\sigma(r)=0$ for all $r$. As a consequence, Definition~\ref{Def-fact-sfnxy} 
takes the form
\begin{equation}
   s^{R'}_{\lambda/\mu}(\x,\y|\a)= \sum_{T\in{\cal T}^{\lambda/\mu}_{R'}} \prod_{(i,j)\in F^{\lambda/\mu}}  \wgt(t_{ij})
\end{equation} 
where $t_{ij}$ is the entry in the box in $i$th row and $j$th column of $T$, with 
\begin{equation}
               \wgt(t_{ij}) = \begin{cases}
							                     x_k+a_{j-i}&\mbox{if $t_{ij}=2k=i_k\in M$};\cr
																	 y_\ell-a_{j-i}&\mbox{if $t_{ij}=(2\ell-1)'=j'_\ell\in N'$}.\cr
							\end{cases}
\end{equation}
This is just the specialisation of the ninth variation to the sixth variation factorial supersymmetric skew Schur functions 
in which $x_{kc}=x_k+a_c$ and $y_{\ell c}=y_\ell-a_c$, and therefore satisfies all the outside decomposition
determinantal identities of Theorem~\ref{The-HGXY}. It is a direct analogue of the factorial skew Schur $Q$-functions that share exactly the same
weighting and satisfy outside decomposition Pfaffian identities~\cite{FK20}.

By way of a final comment, we might point out that one can pass from a ninth variation to what we might call a tenth variation merely by 
replacing $\X=(x_{kc})$ and $\Y=(y_{\ell c})$, with $c=j-i$, by $\X=(x_{k,(i,j)})$ and $\Y=(y_{\ell,(i,j)})$, respectively, in Definition~\ref{Def-sfnXY-R'}, 
just as suggested by Bachmann in writing down his most general non-diagonal form with $c$ replaced by $(i,j)$. However,
such cases do not in general satisfy any of our determinantal identities, and of course display no supersymmetry properties.

\noindent{\bf Acknowledgements}
\label{sec:ack}
The first author (AMF) was supported by a
Discovery Grant from the Natural Sciences and Engineering Research Council of
Canada (NSERC). The second author (RCK) is grateful for the hospitality extended to him 
by Professor Bill Chen at the Center for Applied Mathematics at Tianjin University 
and for the opportunity to pursue this project while visiting him there.
This work was supported by the Canadian Tri-Council Research
Support Fund.

\end{document}